\documentclass[smallextended,referee,envcountsect,]{svjour3}
\smartqed
\usepackage{graphicx}

\usepackage{pgfplots}
\usepackage{epstopdf}%
\usepackage{multirow}%
\usepackage{amssymb,amsfonts}%
\usepackage{mathrsfs}%
\usepackage[title]{appendix}%
\usepackage{xcolor}%
\usepackage{textcomp}%
\usepackage{manyfoot}%
\usepackage{booktabs}%
\usepackage{algorithm}%
\usepackage{algorithmicx}%
\usepackage{algpseudocode}%
\usepackage{listings}%
\usepackage{ragged2e}
\usepackage{fullpage}
\usepackage{changepage}
\usepackage{stfloats}
\usepackage[caption=false]{subfig}
\usepackage{extarrows}
\usepackage{hyperref}
\usepackage{cleveref}
\usepackage{tikz}
\usepackage{pifont}
\usepackage{enumitem}
\usepackage{bm}

\journalname{JOTA}

\renewcommand{\thesubsection}{\textbf{\arabic{section}.\arabic{subsection}}} 

\begin{document}
	
	\title{A General Mixed-Order Primal-Dual Dynamical System with Tikhonov Regularization}

	\author{Hong-lu Li\and Rong Hu\and Xin He\and Yi-bin Xiao}

	\institute{Hong-lu Li \at
		University of Electronic Science and Technology of China \\
		Chengdu, Sichuan 611731, P.R. China \\
		lihonglu$\_$7988@163.com
		\and 
		Rong Hu Corresponding author\at
		Chengdu University of Technology\\
		Chengdu, Sichuan 610059, P.R. China\\
		hrong1130@foxmail.com
		\and
		Xin He \at
		Xihua University \\
		Chengdu, Sichuan 610039, P.R. China \\
		hexinuser@163.com
		\and
		Yi-bin Xiao  Corresponding author  \at
		University of Electronic Science and Technology of China \\
		Chengdu, Sichuan 611731, P.R. China \\
		xiaoyb9999@hotmail.com
	}
	
	\date{Received: date / Accepted: date}
	
	\maketitle
	
	\begin{abstract}
		In a Hilbert space, we propose a class of general mixed-order primal-dual dynamical systems with Tikhonov regularization for a convex optimization problem with linear equality constraints. The proposed dynamical system is characterized by three time-dependent parameters, i.e., general viscous damping, time scaling, and Tikhonov regularization coefficients, which can incorporate as special cases some existing mixed-order primal-dual dynamical systems in the literature. With some appropriate conditions on the parameters, we analyze by constructing suitable Lyapunov functions the asymptotic convergence properties of the proposed dynamical system, where a convergence rate of $\mathcal{O}(\frac{1}{t^2\beta(t)})$ for the objective function error and a convergence rate of $o(\frac{1}{\beta(t)})$ for the primal-dual gap are established. Moreover, we further prove the strong convergence of the trajectory generated by the proposed dynamical system. Finally, we carry out some numerical experiments to illustrate the obtained theoretical results of the proposed dynamical system.
	\end{abstract}
	\keywords{Convex optimization problem with linear equality constraints \and mixed-order primal-dual dynamical system \and Lyapunov function \and asymptotic convergence \and strong convergence}
	\subclass{34D05 \and  37N40 \and 46N10 \and 65B99 \and 90C25}
		
		\noindent\textbf{Communicated by Hong-Kun Xu.}
	
	
	\section{Introduction}\label{sec1}
	\hspace{1.5em}Let $\mathcal{X},\mathcal{Y}$ be two real Hilbert spaces with the inner product $\langle\cdot,\cdot\rangle$ and its induced norm $\Vert \cdot \Vert$. In this paper, we consider a class of linear equality constrained convex optimization problems (convex optimization problems with linear equality constraints) as follows
	\begin{equation}\label{equ1}
		\begin{aligned}
			&\min_{x\in\mathcal{X}}\;\;f(x),\\
			&\;{\rm{s.t.}} \;Ax =b,
		\end{aligned}
	\end{equation}
	where $f:\mathcal{X}\rightarrow\mathbb{R}$ is a continuously differentiable convex function such that $\nabla f$ is Lipschitz continuous,  $A:\mathcal{X}\rightarrow\mathcal{Y}$ is a continuous linear operator, and $b\in\mathcal{Y}$ is given. In what follows, we always consider that the feasible set is nonempty. Due to their important applications in various fields such as machine learning \cite{bib26}, image recovery \cite{bib28}, network optimization \cite{bib29,bib13}, transportation \cite{bib53,bib54}, engineering design \cite{bib55}, healthcare \cite{bib56,bib57}, etc., the linear equality constrained convex optimization problems have been a research hot spot and attracted the attention of many scholars in different fields. 
	
   In the literature, a large number of scholars are devoted to introducing different dynamical systems to study optimization problems since dynamical system is an efficient tool to interpret and design algorithms for optimization problems. The gradient flow dynamical system proposed in \cite{bib58} corresponds to the classic gradient descent algorithm for the unconstrained optimization problem $\min_{x}\; \phi(x)$ with  $\phi $ being a smooth function.
	
	To further improve the convergence rate of algorithms, second-order dynamical systems, which incorporate an inertial term $\ddot{x}(t)$, are proposed to study optimization problems. In 1964, Polyak \cite{bib1} proposed a second-order dynamical system with friction, called heavy ball system, and designed its corresponding heavy ball algorithm, for an unconstrained optimization problem. He further proved that the heavy ball algorithm can attain a faster rate of local convergence than the classic gradient descent method near the minimum of optimization problem. 
	
	In order to understand the Nesterov's accelerated gradient method, which is a classic numerical method for unconstrained optimization problems, Su et al. \cite{bib7} introduced in 2016 the following dynamical system
	\begin{equation*}
		({\rm{AVD}})_{\alpha}\qquad
		\ddot{x}(t)+\frac{\alpha}{t}\dot{x}(t)+\nabla \phi(x(t))=0.
	\end{equation*}
	They showed that the introduced dynamical system  ${\rm{(AVD)_\alpha}}$ with $\alpha=3$ can be viewed as a continuous limit of Nesterov's accelerated gradient algorithm.  Using Lyapunov analysis, they established $\mathcal{O}(1/t^2)$ convergence rate for the objective function error along the trajectory, which matches that of discrete Nesterov's accelerated gradient algorithm. 	  This approach provided a deeper understanding of Nesterov's accelerated gradient algorithm from a dynamical systems perspective, bridging the gap between continuous systems and discrete algorithms. Based on the dynamical system  ${\rm{(AVD)_\alpha}}$, May \cite{bib8} and Attouch et al. \cite{bib9} further proved that the convergence rate of the function value is $o(1/t^2)$ when $\alpha > 3$ and $\mathcal{O}(t^{-2\alpha/3})$ when $\alpha \leq 3$, respectively.

	In the literature, there are many extensions or generalizations of the dynamical system  ${\rm{(AVD)_\alpha}}$ introduced to study unconstrained optimization problems. One important generalization is related to the viscous damping parameter  and the  time scaling parameter in the dynamical system. Cabot et al. \cite{bib59}  introduced a dynamical system featuring a generalized viscous damping coefficient and proved the  strong convergence of the trajectory generated by this dynamical system under some more stringent conditions. Building on this, Attouch et al. \cite{bib12} proposed the following dynamical system in 2019, incorporating  more general viscous damping parameter  $\gamma(t)$ and time scaling parameter $\beta(t)$:
	\begin{equation*}
		\ddot{x}(t)+\gamma(t)\dot{x}(t)+\beta(t)\nabla \phi(x(t))=0.
	\end{equation*}
	The effect of the introduction of more general viscous damping parameter and time scaling parameter is to get a better convergence rate of the objective function error. In \cite{bib12}, Attouch et al. established convergence rate $\mathcal{O}(1/{\beta(t)\Gamma(t)^2})$ for the objective function error along the trajectory, where $\Gamma(t)=e^{\int_{t_0}^{t}\gamma(u)du}\int_{t}^{+\infty}e^{-\int_{t_0}^{u}\gamma(s)ds}du$. Therefore, the convergence rate can be better than existing results in some special cases.  More related results on second-order dynamical systems can be found in \cite{bib3,bib4,bib10,bib11,bib41,bib42,bib43,bib44,bib45,bib46,bib47,bib48,bib50,bib51}. 
	Another generalization of the dynamical system  ${\rm{(AVD)_\alpha}}$ is the introducing of a term which is called Tikhonov regularization term. The Tikhonov regularization technique has been applied to first-order and second-order dynamical systems since it can ensure the strong convergence of the trajectory generated by the dynamical system under appropriate assumptions. Attouch et al. \cite{bib21} proposed the Tikhonov regularization of dynamical system  ${\rm{(AVD)_\alpha}}$ by incorporating  a Tikhonov regularization term $\epsilon(t)x(t)$ as follows
	
	\begin{equation*}
		({\rm{AVD}}){\alpha,\epsilon}\qquad
		\ddot{x}(t)+\frac{\alpha}{t}\dot{x}(t)+\nabla \phi(x(t))+\epsilon(t)x(t)=0.
	\end{equation*}
	Some convergence results are obtained when the Tikhonov regularization parameter $\epsilon(t)$ approaches zero at different rates: (i) when $\epsilon(t)$ approaches zero rapidly, the dynamical system $({\rm{AVD}})_{\alpha,\epsilon}$ has the same fast convergence as the dynamical system $({\rm{AVD}})_{\alpha}$, (ii) when $\epsilon(t)$  approaches zero slowly, the trajectory generated by $({\rm{AVD}})_{\alpha,\epsilon}$ strongly converge to the minimum norm solution of its corresponding unconstrained optimization problem. Following that,  Xu et al. \cite{bib22} considered the Tikhonov regularization of dynamical system  ${\rm{(AVD)_\alpha}}$ with a general time scaling coefficient for unconstrained optimization problem, which can get the strong convergence of its trajectory with better convergence rate under some appropriate assumptions. For more related research studies on Tikhonov regularization of dynamical systems for unconstrained optimization problems, we refer to references \cite{bib18,bib19,bib20,bib35,bib36,bib37,bib38,bib39,bib40,bib49}.

	As we know, transforming it into an unconstrained optimization problem is a crucial and effective method to solve a constrained optimization problem and thus, dynamical systems are also employed to study constrained optimization problems. In this context, we introduce its augmented  Lagrangian function  $\mathcal{L}_{\sigma}:\mathcal{X}\times\mathcal{Y}\rightarrow\mathbb{R}$ as follows
	\begin{equation*}
		\mathcal{L}_{\sigma}(x,\lambda) := f(x)+\langle \lambda,Ax-b\rangle+\frac{\sigma}{2}\Vert Ax-b\Vert^2,
	\end{equation*}
	where $\lambda\in\mathcal{Y}$ is called the Lagrange multiplier or the dual variable and $\sigma\ge 0$ is the penalty parameter. The linear equality constrained convex optimization problem \eqref{equ1} corresponds to the following saddle point problem
	\begin{equation}\label{equ2-1}
		\min_{x\in \mathcal{X}}\max_{\lambda\in \mathcal{Y}}\mathcal{L}_{\sigma}(x,\lambda).
	\end{equation}	
	In this paper, we assume that the linear equality constrained convex optimization problem \eqref{equ1} satisfies the Slater condition. Then $x^*$ is optimal for the problem \eqref{equ1} and $\lambda^*$ is a corresponding Lagrange multiplier if and only if $(x^*, \lambda^*)$ is a saddle point of $\mathcal{L}_{\sigma}$, i.e.,
	\begin{equation*}\label{equ2}
		\mathcal{L}_{\sigma}(x^*,\lambda)\leq\mathcal{L}_{\sigma}(x^*,\lambda^*)\leq\mathcal{L}_{\sigma}(x,\lambda^*),\quad \forall(x,\lambda)\in\mathcal{X}\times\mathcal{Y}.
	\end{equation*}
	We denote by $\Omega$ the set of saddle points of the problem \eqref{equ2-1} (primal-dual optimal solution of the problem \eqref{equ1}), which is assumed to be nonempty in the following. By using the above saddle point problem, many scholars have proposed several second-order dynamical systems for solving the linear equality constrained convex optimization problem \eqref{equ1}. Zeng et al.  \cite{bib13} first introduced and studied the following second-order dynamical system for solving the problem \eqref{equ1}
	\begin{equation*}
		\begin{cases}
			\ddot{x}(t)+\frac{\alpha}{t}\dot{x}(t)=-\nabla f(x(t))-A^T(\lambda(t)+\beta t\dot{\lambda}(t))-A^T(Ax(t)-b),\\
			\ddot{\lambda}(t)+\frac{\alpha}{t}\dot{\lambda}(t)=A(x(t)+\beta t\dot{x}(t))-b.
		\end{cases}
	\end{equation*}
	They further proved that $\mathcal{L}(x(t),\lambda^*)-\mathcal{L}(x^*,\lambda^*)=\mathcal{O}(t^{-\frac{2\min\left\{3,\alpha\right\}}{3}})$ when $\alpha>0$ and $\beta=\frac{3}{2\min\left\{3,\alpha\right\}}$. He et al.  \cite{bib14} and Attouch et al.  \cite{bib15} each proposed a generalized second-order dynamical system with time-dependent damping coefficients for solving a separable convex optimization problem with linear equality constraints. Both works analyzed the fast convergence of the primal-dual gap and the feasibility violation along their respective  trajectory. Based on their obtained results, Bo\c{t} and Nguyen \cite{bib16} further proved that the primal-dual trajectory of the introduced second-order dynamical system weakly converges to the primal-dual solution of the saddle point problem corresponding to the separable convex optimization problem. It should be noted that the aforementioned second-order dynamical systems for constrained optimization problems involve the second-order terms of both primal and dual variables. In contrast, He et al.  \cite{HeTAC2022} first proposed a ``second-order" + ``first-order" primal-dual dynamical system for the  linear equality constrained convex optimization problem \eqref{equ1} with constant damping. Later, He et al. \cite{bib17} extended the results of \cite{HeTAC2022} to the case of viscous damping $\frac{\alpha}{t}$, and proposed the dynamical system for solving the problem \eqref{equ1} as follows 
	\begin{align}
		\begin{cases}
			\ddot{x}(t)+\frac{\alpha}{t}\dot{x}(t)&=-\beta(t)(\nabla f(x(t))-A^T\lambda(t)),\\
			\dot{\lambda}(t)&=t\beta(t)(A(x(t)+\frac{t}{\alpha-1}\dot{x}(t))-b),
		\end{cases} \nonumber
	\end{align}
	where $\alpha>1$ and $\beta(t)$ is a positive time scaling coefficient. From the viewpoint of the numerical computation, a pure second-order dynamical system is more challenging than a mixed-order (``second-order" + ``first-order") dynamical system. Nevertheless, by using the introduced mixed-order dynamical system, He et al. \cite{bib17} achieved the same convergence rate as the case with a pure second-order dynamical system.
	Based on this, Zhu et al. \cite{bib23} recently considered the Tikhonov regularization of a mixed-order primal-dual dynamical system for the linear equality constrained convex optimization problem \eqref{equ1}, and there, the strong convergence results on the trajectory of the dynamical system are obtained.
	
	Inspired by the above studies on the dynamical systems for optimization problems, we are devoted to proposing a general dynamical system for the linear equality constrained convex optimization problem \eqref{equ1} so as to improve the convergence properties of such dynamical systems. To this end, we consider a general mixed-order primal-dual dynamical system with Tikhonov regularization for the problem \eqref{equ1} as follows
	
	\begin{align}
		&	\begin{cases}
			\begin{aligned}
				&\ddot{x}(t) + \gamma(t) \dot{x}(t) + \beta(t) (\nabla f(x(t)) + A^T \lambda(t) + \sigma A^T (Ax(t) - b) + \epsilon(t) x(t)) = 0, \\
				&\dot{\lambda}(t) - t \beta(t) (A(x(t) + \theta t \dot{x}(t)) - b) = 0,
			\end{aligned}
		\end{cases} \label{equ5} \\
		\intertext{\vspace{4ex} i.e.,\vspace{-8ex}} \nonumber \\
		&\begin{cases}
			\begin{aligned}
				&\ddot{x}(t) + \gamma(t) \dot{x}(t) + \beta(t) (\nabla_x \mathcal{L}_{\sigma}(x(t), \lambda(t)) + \epsilon(t) x(t)) = 0,\\
				&\dot{\lambda}(t) - t \beta(t) \nabla_{\lambda} \mathcal{L}_{\sigma}(x(t) + \theta t \dot{x}(t), \lambda(t)) = 0,
			\end{aligned}
		\end{cases} \nonumber
	\end{align}
	where $\theta>0$, $t\ge t_0>0$, $\gamma,\beta:\left[t_0,+\infty\right)\rightarrow\left(0,+\infty\right)$ are two continuous functions, and $\epsilon:\left[t_0,+\infty\right)\rightarrow\mathbb{R}_+$ is a $\mathcal{C}^1$ non-increasing function  such that $\lim\limits_{t\rightarrow+\infty}\epsilon(t)=0$. Obviously, the introduced dynamical system \eqref{equ5} for the linear equality constrained convex optimization problem \eqref{equ1} includes as special cases the dynamical systems for unconstrained and constrained optimization problems. For instance, see \cite{bib7,bib21,bib22,bib23}.
	
	This paper is organized as follows. In Section \ref{sec2}, by constructing appropriate Lyapunov functions, we analyze the asymptotic convergence properties of the dynamical system \eqref{equ5} for the linear equality constrained convex optimization problem \eqref{equ1}, based on which, we discuss in Section \ref{sec3} the strong convergence of the trajectory generated by the dynamical system \eqref{equ5} when the regularization parameter $\epsilon(t)$ approaches zero at a suitable rate. In Section \ref{sec4}, we carry out some numerical experiments to illustrate the obtained theoretical results. Finally, Section \ref{sec5} concludes the paper by summarizing the key findings and contributions of this research.

	\section{Asymptotic Analysis}\label{sec2}
	In this section, we focus on the asymptotic convergence properties of the dynamical system \eqref{equ5}, where the existence and uniqueness of the global solution for its corresponding Cauchy problem can be easily proven by applying the Cauchy-Lipschitz-Picard theorem and we refer the readers to Appendix \ref{Appendix A} for a detailed proof. We first prove some convergence results for the dynamical system \eqref{equ5} by considering two cases where the Tikhonov regularization parameter $\epsilon(t)$ approaches zero at different rates. Subsequently, we present specific examples corresponding to these two cases.
	
	\subsection{\bf Case $\bm{\int_{t_0}^{+\infty}t\beta(t)\epsilon(t) dt<+\infty}$}
	In this subsection, we analyze the asymptotic behavior of the dynamical system \eqref{equ5} with $\beta(t)$ satisfying $\int_{t_0}^{+\infty}t\beta(t)\epsilon(t) dt<+\infty$, which indicates that the Tikhonov regularization parameter $\epsilon(t)$ converges to 0 at a relatively fast rate. To proceed, based on the introduced dynamical system \eqref{equ5}, we first define a function $\mathcal{G}:\left[t_0,+\infty\right)\rightarrow\left[0,+\infty\right)$ as follows: 
	\begin{equation}\label{equ32}
		\begin{aligned}
			\mathcal{G}(t) :=\;& b(t)^2\beta(t)(\mathcal{L}_{\sigma}(x(t),\lambda^*)-\mathcal{L}_{\sigma}(x^*,\lambda^*)+\frac{\epsilon(t)}{2}\Vert x(t)\Vert^2)+\frac{c(t)}{2}\Vert x(t)-x^*\Vert^2\\
			& +\frac{1}{2}\Vert \eta(x(t)-x^*)+b(t)\dot{x}(t)\Vert^2+\frac{1}{2}\Vert\lambda(t)-\lambda^*\Vert^2,
		\end{aligned}
	\end{equation}
	where $b(t)$ := $\sqrt{\theta}t$, $c(t)$ := $t\gamma(t)-\frac{1+\theta}{\theta}$ and $ \eta$ := $\frac{1}{\sqrt{\theta}}$. The function $\mathcal{G}$ is called a time-varying Lyapunov function for the dynamical system \eqref{equ5} since it is non-negative and bounded under some assumptions which will be specified later.

	\begin{lemma}\label{lem2.1}
		
		Let $\mathcal{X},\mathcal{Y}$ be two real Hilbert spaces, $f:\mathcal{X}\rightarrow\mathbb{R}$ be a continuously differentiable convex function such that $\nabla f$ is Lipschitz continuous,  $A:\mathcal{X}\rightarrow\mathcal{Y}$ be a continuous linear operator, and $b\in\mathcal{Y}$ be a given point in $\mathcal{Y}$, and  let $\theta>0$ and $t_0>0$ be two constants. Assume that $\gamma,\beta,\epsilon:\left[t_0,+\infty\right)\rightarrow\left[0,+\infty\right)$ are three differentiable functions. Then, for any trajectory (global solution) $(x(t),\lambda(t))$ of the dynamical system $\rm{\eqref{equ5}}$ and any primal-dual optimal solution $(x^*,\lambda^*)\in\Omega$ of the linear equality constrained optimization problem $\rm{\eqref{equ1}}$,
		\begin{equation}\label{equ33}
			\begin{aligned}
				\dot{\mathcal{G}}(t)\leq\;&((2\theta-1)t\beta(t)+\theta t^2\dot{\beta}(t))(\mathcal{L}_{\sigma}(x(t),\lambda^*)-\mathcal{L}_{\sigma}(x^*,\lambda^*))+\frac{t\beta(t)\epsilon(t)}{2}\Vert x^*\Vert^2\\
				&+\frac{1}{2}((2\theta-1)t\beta(t)\epsilon(t)+\theta t^2\dot{\beta}(t)\epsilon(t)+\theta t^2\beta(t)\dot{\epsilon}(t))\Vert x(t)\Vert^2-\frac{\sigma t\beta(t)}{2}\Vert Ax(t)-b\Vert^2\\
				&+(\theta+1-\theta t\gamma(t))t\Vert\dot{x}(t)\Vert^2+\frac{1}{2}(\gamma(t)+t\dot{\gamma}(t)-t\beta(t)\epsilon(t))\Vert x(t)-x^*\Vert^2.
			\end{aligned}
		\end{equation}
	\end{lemma}
	{\it Proof} 
	By calculating the derivative of \eqref{equ32}, we obtain
	\begin{equation}\label{equ45}
		\begin{split}
			\dot{\mathcal{G}}(t)=& (2b(t)\dot{b}(t)\beta(t)+b(t)^2\dot{\beta}(t))(\mathcal{L}_{\sigma}(x(t),\lambda^*)-\mathcal{L}_{\sigma}(x^*,\lambda^*))+b(t)^2\beta(t)\langle\nabla_x\mathcal{L}_{\sigma}(x(t),\lambda^*),\dot{x}(t)\rangle\\
			&+\frac{1}{2}((2b(t)\dot{b}(t)\beta(t)+b(t)^2\dot{\beta}(t))\epsilon(t)+b(t)^2\beta(t)\dot{\epsilon}(t))\Vert x(t)\Vert^2+b(t)^2\beta(t)\epsilon(t)\langle x(t),\dot{x}(t)\rangle\\
			&+\langle \eta(x(t)-x^*)+b(t)\dot{x}(t),b(t)\ddot{x}(t)\rangle+(\eta^2+\eta\dot{b}(t)+c(t))\langle x(t)-x^*,\dot{x}(t)\rangle\\
			&+(\eta b(t)+b(t)\dot{b}(t))\Vert\dot{x}(t)\Vert^2+\frac{\dot{c}(t)}{2}\Vert x(t)-x^*\Vert^2+\langle\lambda(t)-\lambda^*,\dot{\lambda}(t)\rangle.
		\end{split}
	\end{equation}
	Meanwhile, by using the definition of the dynamical system \eqref{equ5}, it follows that
	\begin{align*}
		&\langle \eta(x(t)-x^*)+b(t)\dot{x}(t),b(t)\ddot{x}(t)\rangle\\
		=& \langle \eta b(t)(x(t)-x^*)+b(t)^2\dot{x}(t),-\gamma(t)\dot{x}(t)-\beta(t)(\nabla_x\mathcal{L}_{\sigma}(x(t),\lambda(t))+\epsilon(t)x(t))\rangle\\
		=&-\eta b(t)\gamma(t)\langle x(t)-x^*,\dot{x}(t)\rangle-\eta b(t)\beta(t)\epsilon(t)\langle x(t)-x^*,x(t)\rangle-b(t)^2\gamma(t)\Vert\dot{x}(t)\Vert^2\\
		&-b(t)^2\beta(t)\epsilon(t)\langle x(t),\dot{x}(t)\rangle-\eta b(t)\beta(t)\langle \nabla_x\mathcal{L}_{\sigma}(x(t),\lambda(t)), x(t)-x^*\rangle\\
		&-b(t)^2\beta(t)\langle \nabla_x\mathcal{L}_{\sigma}(x(t),\lambda(t)),\dot{x}(t)\rangle,
	\end{align*}
	and
	\begin{equation*}
		\langle\lambda(t)-\lambda^*,\dot{\lambda}(t)\rangle=t\beta(t)\langle A^T(\lambda(t)-\lambda^*),x(t)-x^*+\theta t\dot{x}(t)\rangle.
	\end{equation*}
	The above two equations are substituted into \eqref{equ45}, thereby yielding the following equation
	\begin{equation}\label{equ9}
		\begin{aligned}
			\dot{\mathcal{G}}(t) =\;&(2b(t)\dot{b}(t)\beta(t)+b(t)^2\dot{\beta}(t))(\mathcal{L}_{\sigma}(x(t),\lambda^*)-\mathcal{L}_{\sigma}(x^*,\lambda^*))-\eta b(t)\beta(t)\epsilon(t)\langle x(t)-x^*,x(t)\rangle\\
			&+\frac{\dot{c}(t)}{2}\Vert x(t)-x^*\Vert^2-\eta b(t)\beta(t)\langle\nabla_x\mathcal{L}_{\sigma}(x(t),\lambda^*),x(t)-x^*\rangle\\
			&+(\eta^2+\eta\dot{b}(t)-\eta b(t)\gamma(t)+c(t))\langle x(t)-x^*,\dot{x}(t)\rangle\\
			&+(b(t)\dot{b}(t)\beta(t)\epsilon(t)+\frac{1}{2}(b(t)^2\dot{\beta}(t)\epsilon(t)+b(t)^2\beta(t)\dot{\epsilon}(t)))\Vert x(t)\Vert^2\\
			&+(\theta t^2\beta(t)-b(t)^2\beta(t))\langle\lambda(t)-\lambda^*,A\dot{x}(t)\rangle+(t\beta(t)-\eta b(t)\beta(t))\langle\lambda(t)-\lambda^*,Ax(t)-b\rangle\\
			&+(\eta b(t)+b(t)\dot{b}(t)-b(t)^2\gamma(t))\Vert\dot{x}\Vert^2.
		\end{aligned}
	\end{equation}
	Since $f(x)+\frac{\epsilon(t)}{2}\Vert x\Vert^2$ is an $\epsilon(t)$-strongly convex function, it follows that
	\vspace{-5pt}
	\begin{equation*}
		f(x^*)+\frac{\epsilon(t)}{2}\Vert x^*\Vert^2-f(x(t))-\frac{\epsilon(t)}{2}\Vert x(t)\Vert^2\\
		\ge\langle\nabla f(x(t))+\epsilon(t)x(t),x^*-x(t)\rangle+\frac{\epsilon(t)}{2}\Vert x(t)-x^*\Vert^2.
	\end{equation*}
	As a consequence,
	\begin{equation}\label{equ10}
		\begin{aligned}
			\langle\nabla_x\mathcal{L}_{\sigma}(x(t),\lambda^*),x(t)-x^*\rangle =\;&\langle \nabla f(x(t))+A^T\lambda^*+\sigma A^T(Ax(t)-b),x(t)-x^*\rangle\\
			=\;&\langle \nabla f(x(t)),x(t)-x^*\rangle+\langle\lambda^*,A(x(t))-b\rangle+\sigma\Vert Ax(t)-b\Vert^2\\
			\ge \;& (\mathcal{L}_{\sigma}(x(t),\lambda^*)-\mathcal{L}_{\sigma}(x^*,\lambda^*))+\frac{\epsilon(t)}{2}(\Vert x(t)-x^*\Vert^2+\Vert x(t)\Vert^2-\Vert x^*\Vert^2)\\
			\;&-\epsilon(t)\langle x(t),x(t)-x^*\rangle+\frac{\sigma}{2}\Vert Ax(t)-b\Vert^2.
		\end{aligned}
	\end{equation}
	By combining \eqref{equ10}, we can get from \eqref{equ9} that
	\begin{equation}\label{equ11}
		\begin{aligned}
			\dot{\mathcal{G}}(t) \leq\;&(2b(t)\dot{b}(t)\beta(t)+b(t)^2\dot{\beta}(t)-\eta b(t)\beta(t))(\mathcal{L}_{\sigma}(x(t),\lambda^*)-\mathcal{L}_{\sigma}(x^*,\lambda^*))\\
			&+\frac{1}{2}(2b(t)\dot{b}(t)\beta(t)\epsilon(t)+b(t)^2\dot{\beta}(t)\epsilon(t)+b(t)^2\beta(t)\dot{\epsilon}(t)-\eta b(t)\beta(t)\epsilon(t))\Vert x(t)\Vert^2\\
			&+\frac{1}{2}(\dot{c}(t)-\eta b(t)\beta(t)\epsilon(t))\Vert x(t)-x^*\Vert^2+\frac{\eta b(t)\beta(t)\epsilon(t)}{2}\Vert x^*\Vert^2\\
			&+ (\eta ^2+\eta \dot{b}(t)-\eta b(t)\gamma(t)+c(t))\langle x(t)-x^*,\dot{x}(t)\rangle\\
			&+(\eta b(t)+b(t)\dot{b}(t)-b(t)^2\gamma(t))\Vert\dot{x}(t)\Vert^2-\frac{\sigma\eta b(t)\beta(t)}{2}\Vert Ax(t)-b\Vert^2\\
			&+(\theta t^2\beta(t)-b(t)^2\beta(t))\langle\lambda(t)-\lambda^*,A\dot{x}(t)\rangle+ (t\beta(t)-\eta b(t)\beta(t))\langle\lambda(t)-\lambda^*,Ax(t)-b\rangle.
		\end{aligned}
	\end{equation}
	Furthermore, based on the specified formulations of $\eta$, $b(\cdot)$, and $c(\cdot)$, it follows that
	\begin{equation*}\label{equ12}
		\begin{cases}
			\eta^2+\eta\dot{b}(t)-\eta b(t)\gamma(t)+c(t)=0,\\[3pt]
			\theta t^2\beta(t)-b(t)^2\beta(t)=0,\\[3pt]
			t\beta(t)-\eta b(t)\beta(t)=0.
		\end{cases}
	\end{equation*}
	Therefore, this together with \eqref{equ11} indicates that  \eqref{equ33} holds.
	This completes the proof of Lemma \ref{lem2.1}.
	\qed
	\vspace{-2pt}
	\begin{theorem}\label{thm2.1}
		Let all hypotheses in Lemma \ref{lem2.1} hold. Moreover, assume further that  $\epsilon$ is a $\mathcal{C}^1$ and non-increasing function such that $\int_{t_0}^{+\infty}t\beta(t)\epsilon(t) dt<+\infty$, $\gamma$ is a $\mathcal{C}^1$ function, $\beta$ is a $\mathcal{C}^2$ and non-negative function satisfying $\lim\limits_{t\rightarrow+\infty}t^2\beta(t)=+\infty$, and
		\begin{eqnarray}
			&&(2\theta-1)\beta(t)+\theta t\dot{\beta}(t)\leq 0,\;(\forall t\ge t_0)\label{equ6}\\[5pt]
			&&\gamma(t)+t\dot{\gamma}(t)\leq t\beta(t)\epsilon(t),\;(\forall t\ge t_0)\label{equ7}\\[5pt]
			&&\theta t \gamma(t)-\theta-1\ge 0.\;(\forall t\ge t_0)\label{equ8}\vspace{-8pt}
		\end{eqnarray}
		Then, for any trajectory $(x(t),\lambda(t))$ of the dynamical system $\rm{\eqref{equ5}}$ and any primal-dual optimal solution $(x^*,\lambda^*)\in\Omega$ of the problem $\rm{\eqref{equ1}}$, the following conclusions hold:
		\begin{flalign*}
		& \begin{array}{l@{\quad}l}
			\textbf{(i)\;(Boundedness of Trajectory)} & (x(\cdot),\lambda(\cdot))\text{ is bounded on } [t_0, +\infty);\\[5pt]
			\textbf{(ii)\;(Pointwise Estimates)} & 
			\left\{
			\begin{array}{l}
				\Vert \dot{x}(t)\Vert=\mathcal{O}(\frac{1}{t}),\\[5pt]
				\mathcal{L}_{\sigma}(x(t),\lambda^*)-\mathcal{L}_{\sigma}(x^*,\lambda^*)=\mathcal{O}(\frac{1}{t^2\beta(t)}),\\[5pt]
				\Vert  \nabla f(x(t))-\nabla f(x^*)\Vert=\mathcal{O}(\frac{1}{t\sqrt{\beta(t)}}),\\[5pt]
				\Vert Ax(t)-b\Vert=\mathcal{O}(\frac{1}{t^2\beta(t)}),\\[5pt]
				\vert f(x(t))-f(x^*)\vert =\mathcal{O}(\frac{1}{t^2\beta(t)});
			\end{array}
			\right. \vspace{5pt}
			\\
						\textbf{(iii)\;(Integral Estimates)} & 
			\begin{array}{l}	
				\int_{t_0}^{+\infty}t\beta(t)\Vert Ax(t)-b\Vert^2 dt<+\infty;\\[5pt]
			\end{array}
		\end{array} &
		\end{flalign*}
		\indent In particular, if $(2\theta-1)\beta(t)+\theta t\dot{\beta}(t)< 0,\;\forall t\ge t_0$, then 
		\begin{equation*}
			\begin{cases}
				\int_{t_0}^{+\infty} ((1-2\theta)\beta(t)-\theta t\dot{\beta}(t))t	(\mathcal{L}_{\sigma}(x(t),\lambda^*)-\mathcal{L}_{\sigma}(x^*,\lambda^*))dt<+\infty,\\[5pt]
				\int_{t_0}^{+\infty} ((1-2\theta)\beta(t)-\theta t\dot{\beta}(t))t \Vert  \nabla f(x(t))-\nabla f(x^*)\Vert^2 dt<+\infty;
			\end{cases}
		\end{equation*}
	and if $\theta t\gamma(t)-\theta-1>0,\;\forall t\ge t_0$, then $	\int_{t_0}^{+\infty}(\theta t\gamma(t)-\theta-1)t\Vert \dot{x}(t)\Vert dt<+\infty$.
	\end{theorem}
	{\it Proof} 
	According to Lemma \ref{lem2.1}, integrating \eqref{equ33} from $t_0$ to $t$ yields
	\begin{equation}\label{equ13}
		\begin{aligned}
			\mathcal{G}(t)&-\int_{t_0}^{t}((2\theta -1)\tau\beta(\tau)+\theta\tau^2\dot{\beta}(\tau))(\mathcal{L}_{\sigma}(x(\tau),\lambda^*)-\mathcal{L}_{\sigma}(x^*,\lambda^*))  d\tau\\ 
			&-\frac{1}{2}\int_{t_0}^{t}((2\theta-1)\tau\beta(\tau)\epsilon(\tau)+\theta\tau^2\dot{\beta}(\tau)\epsilon(\tau)+\theta\tau^2\beta(\tau)\dot{\epsilon}(\tau))\Vert x(\tau)\Vert^2d\tau\\
			&-\frac{1}{2}\int_{t_0}^{t}(\gamma(\tau)+\tau\dot{\gamma}(\tau)-\tau\beta(\tau)\epsilon(\tau))\Vert x(\tau)-x^*\Vert^2 d\tau\\
			&+\frac{1}{2}\int_{t_0}^{t}\sigma\tau\beta(\tau)\Vert Ax(\tau)-b\Vert^2 d\tau-\int_{t_0}^{t}(\theta+1-\theta\tau\gamma(\tau))\tau\Vert\dot{x}(\tau)\Vert^2 d\tau\\
			\leq&\;\mathcal{G}(t_0) +\int_{t_0}^{t}\frac{\tau\beta(\tau)\epsilon(\tau)}{2}\Vert x^*\Vert^2 d\tau.
		\end{aligned}
	\end{equation}
	Based on the assumptions in Theorem \ref{thm2.1} , it is straightforward to derive the following inequalities
	\begin{equation*}\label{eq-15}
		\begin{cases}
			(2\theta -1)t\beta(t)+\theta t^2\dot{\beta}(t)\leq  0,\\
			\gamma(t)+t\dot{\gamma}(t)-t\beta(t)\epsilon(t)\leq 0,\\
			(\theta +1 -\theta  t\gamma(t))t\leq 0,\\
			(2\theta-1)t\beta(t)\epsilon(t)+\theta t^2\dot{\beta}(t)\epsilon(t)+\theta t^2\beta(t)\dot{\epsilon}(t)\leq 0,\\
			\sigma t\beta(t)\ge 0,
		\end{cases}
	\end{equation*}

  \noindent where the first to third inequalities hold from \eqref{equ6}-\eqref{equ8}, the fourth inequality is satisfied since 
	$\epsilon(t)$ is a non-negative and non-increasing function, and $\beta(t)$ satisfies \eqref{equ6}, and the last inequality holds due to the non-negativity of $\sigma$ and $\beta(t)$. This together with  the assumption  $\int_{t_0}^{+\infty}t\beta(t)\epsilon(t) dt<+\infty$  and \eqref{equ13} implies that 
	$\mathcal{G}(t)$ is bounded on $\left[t_0, +\infty\right)$ (i.e., $\mathcal{G}(t)$ is a Lyapunov function) and $\int_{t_0}^{+\infty}t\beta(t)\Vert Ax(t)-b\Vert^2 dt<+\infty$. In particular, if $\theta t\gamma(t)-\theta-1>0,\;\forall t\ge t_0$, it follows from \eqref{equ13} that
	\begin{equation*}
		\int_{t_0}^{+\infty}(\theta  t\gamma(t)-\theta-1)t\Vert\dot{x}(t)\Vert^2dt<+\infty.
	\end{equation*}
   Furthermore, if $(2\theta-1)\beta(t)+\theta t\dot{\beta}(t)< 0,\;\forall t\ge t_0$, we obtain from \eqref{equ13} that
   \begin{equation}\label{equ99}
   	\int_{t_0}^{+\infty} ((1-2\theta)\beta(t)-\theta t\dot{\beta}(t))t	(\mathcal{L}_{\sigma}(x(t),\lambda^*)-\mathcal{L}_{\sigma}(x^*,\lambda^*))dt<+\infty.
   \end{equation}

	\noindent Moreover, from the definition of $\mathcal{G}(t)$,  there exist two positive constants $C$ and $\widetilde{C}$ such that for all $t \in\left[t_0,+\infty\right)$,
	\begin{equation}\label{eq-2}
		\Vert\lambda(t)-\lambda^*\Vert\leq C\quad\text{and}\quad	\Vert (x(t)-x^*)+\theta t\dot{x}(t)\Vert\leq\sqrt{\theta}\widetilde{C},
	\end{equation}
	which imply that \(\lambda(t)\) and \(x(t)\) are bounded on $\left[t_0, +\infty\right)$, respectively, where we use  Lemma \ref{lemB.3} in Appendix B to get the boundedness of $x(t)$. Therefore, the trajectory $(x(\cdot),\lambda(\cdot))\text{ is bounded on } [t_0, +\infty)$.
	
	Since  $\theta t\Vert\dot{x}(t)\Vert\leq\Vert (x(t)-x^*)+\theta t\dot{x}(t)\Vert+\Vert  (x(t)-x^*)\Vert $, it is conclude by the boundedness of $x(t)$ and \eqref{eq-2} that $\dot{x}(t)$ is bounded and  $\Vert\dot{x}(t)\Vert=\mathcal{O}(\frac{1}{t})$. According to the definition of $\mathcal{G}(t)$, it follows that 
	\begin{equation*}
		t^2\beta(t)(\mathcal{L}_{\sigma}(x(t),\lambda^*)-\mathcal{L}_{\sigma}(x^*,\lambda^*))\leq \mathcal{G}(t), 
	\end{equation*}
	which together with the boundedness of $\mathcal{G}(t)$ and the assumption  $\lim\limits_{t\rightarrow+\infty}t^2\beta(t)=+\infty$ that 
	\begin{eqnarray}
		&&\mathcal{L}_{\sigma}(x(t),\lambda^*)-\mathcal{L}_{\sigma}(x^*,\lambda^*)=\mathcal{O}(\frac{1}{t^2\beta(t)}).\label{equ14}
	\end{eqnarray}
	
	\noindent Since $f$ is convex and $\nabla f$ is Lipschitz continuous, there exists a constant $L > 0$ such that
	\begin{equation*}
		f(x(t)) \ge f(x^*) +\langle \nabla f(x^*), x(t)-x^*\rangle+\frac{1}{2L}\Vert \nabla f(x(t))-\nabla f(x^*)\Vert^2.
	\end{equation*}
	Thus, by the KKT condition at $(x^*, y^*)$ for the problem \eqref{equ1}, we can get that 
	\begin{equation}\label{equ38}
		\begin{aligned}
			\mathcal{L}_{\sigma}(x(t),\lambda^*)-\mathcal{L}_{\sigma}(x^*,\lambda^*)&=f(x(t))-f(x^*)+\langle \lambda^*, Ax(t)-b\rangle+\frac{\sigma}{2}\Vert Ax(t)-b\Vert^2\\
			&\ge f(x(t))-f(x^*)+\langle \lambda^*, Ax(t)-b\rangle\\
			&\ge \langle \nabla f(x^*), x(t)-x^*\rangle+\frac{1}{2L}\Vert \nabla f(x(t))-\nabla f(x^*)\Vert^2+\langle \lambda^*, A(x(t)-x^*)\rangle\\
			& =\frac{1}{2L}\Vert \nabla f(x(t))-\nabla f(x^*)\Vert^2,
		\end{aligned}
	\end{equation}
	which together with \eqref{equ14} implies that 
	\vspace{-5pt}
	\begin{equation*}
		\Vert  \nabla f(x(t))-\nabla f(x^*)\Vert=\mathcal{O}(\frac{1}{t\sqrt{\beta(t)}}).
	\end{equation*}
  In particular, if $(2\theta-1)\beta(t)+\theta t\dot{\beta}(t)< 0,\;\forall t\ge t_0$, then combining this with \eqref{equ13} and \eqref{equ99}, we obtain that
  \begin{equation*}
  		\int_{t_0}^{+\infty} ((1-2\theta)\beta(t)-\theta t\dot{\beta}(t))t \Vert  \nabla f(x(t))-\nabla f(x^*)\Vert^2 dt<+\infty.
  \end{equation*}
	Moreover, it follows from the definition of the dynamical system \eqref{equ5} that
	\begin{equation*}
		\begin{aligned}
			\lambda(t)-\lambda(t_0)=&\int_{t_0}^{t}\dot{\lambda}(\tau)d\tau
			=\int_{t_0}^{t}\tau\beta(\tau)(A(x(\tau)+\theta\tau\dot{x}(\tau))-b)d\tau\\
			=&\int_{t_0}^{t}\tau\beta(\tau)(Ax(\tau)-b)d\tau+\int_{t_0}^{t}\tau\beta(\tau)\theta\tau A\dot{x}(\tau)d\tau\\
			=&\int_{t_0}^{t}\tau\beta(\tau)(Ax(\tau)-b)d\tau+\int_{t_0}^{t}\theta{\tau}^2\beta(\tau) d(Ax(\tau)-b)\\
			=&\theta t^2\beta(t)(Ax(t)-b)-\theta {t_0}^2\beta(t_0)(Ax(t_0)-b)\\
			&+\int_{t_0}^{t}((1-2\theta)\tau\beta(\tau)-\theta\tau^2\dot{\beta}(\tau))(Ax(\tau)-b)d\tau.
		\end{aligned}
	\end{equation*}
	This together with the boundedness of $\lambda(t)$ on $\left[t_0, +\infty\right)$ implies that for all $ t\ge t_0$,
	\begin{equation}\label{equ15}
		\Vert t^2\beta(t)(Ax(t)-b)+\int_{t_0}^{t}(\frac{1-2\theta}{\theta}\tau\beta(\tau)-\tau^2\dot{\beta}(\tau))(Ax(\tau)-b)d\tau \Vert<+\infty.
	\end{equation}
	\noindent Applying Lemma \ref{lemB.1} with $g(t)=t^2\beta(t)(Ax(t)-b)$ and $a(t) = \frac{(1-2\theta)\beta(t)-\theta t\dot{\beta}(t)}{\theta t\beta(t)}$ , we obtain by $\eqref{equ15}$ that
	\begin{equation}\label{eq-3}
		\sup_{t\ge t_0} \Vert t^2\beta(t)(Ax(t)-b)\Vert< +\infty,\quad  i.e.,\, \Vert Ax(t)-b\Vert=\mathcal{O}(\frac{1}{t^2\beta(t)}).
	\end{equation}
	Since 
	\begin{equation*}
		\mathcal{L}_{\sigma}(x(t),\lambda^*)-\mathcal{L}_{\sigma}(x^*,\lambda^*)=f(x(t))-f(x^*)+\langle \lambda^*,Ax(t)-b\rangle+\frac{\sigma}{2}\Vert Ax(t)-b\Vert^2,
	\end{equation*}
	it follows from \eqref{equ14} and \eqref{eq-3} that 
	\begin{equation*}
		\vert f(x(t))-f(x^*)\vert =\mathcal{O}(\frac{1}{t^2\beta(t)}).
	\end{equation*}
	This completes the proof of Theorem \ref{thm2.1}.
	\qed
	\begin{remark}
		It is easy to verify that both the dynamical system $({\rm{AVD}})_{\alpha}$ proposed by Su et al. in \cite{bib7} for unconstrained optimization problems and the Tikhonov regularized second-order plus first-order primal-dual dynamical system proposed by Zhu et al. in \cite{bib23} for  linear equality constrained convex optimization problems are particular cases of the dynamical system \eqref{equ5}. Moreover, by setting the parameters $\gamma(t)=\frac{\alpha}{t}$, $\beta(t)=1$ and $\theta=\frac{1}{\alpha-1}$, we can obtain the same conclusion as Theorem 3.2 in \cite{bib23}.
	  In a related work, He et al. \cite{bib17} obtained fast convergence rates for the objective value error and constraint violation under the weaker condition $\int_{t_0}^{+\infty}t\epsilon(t)dt<+\infty$. These convergence rates match those given in Theorem \ref{thm2.1}. However, the result in Theorem \ref{thm2.1} is not a direct consequence of \cite[Theorem 2]{bib17}, since the boundedness of the trajectory $(x(t))_{t\ge t_0}$ cannot be guaranteed ensured without prior analysis. Therefore, at this stage, the condition $\int_{t_0}^{+\infty}t\beta(t)\epsilon(t)dt<+\infty$ appears to be indispensable in order to establish the conclusions of Theorem 2.1. This observation also suggests potential directions for future research.
	\end{remark}
	
	\subsection{\bf Case $\bm{\int_{t_0}^{+\infty}\frac{\beta(t)\epsilon(t)}{t}  dt<+\infty}$}\label{subsec2}
	In this subsection, we consider the case  with the Tikhonov regularization parameter $\epsilon(t)$ converging  to 0 at a relatively slower rate (i.e., satisfies $\int_{t_0}^{+\infty}\frac{\beta(t)\epsilon(t)}{t} dt<+\infty$).   We redefine the following  Lyapunov function $\widetilde{\mathcal{G}}\colon \left[t_0,+\infty\right)\rightarrow\left[0,+\infty\right)$ to analyze the asymptotic convergence properties of the dynamical system $\eqref{equ5}$.
	\begin{equation}\label{equ36}
		\begin{aligned}
			\widetilde{\mathcal{G}}(t) :=\;& \beta(t)(\mathcal{L}_{\sigma}(x(t),\lambda^*)-\mathcal{L}_{\sigma}(x^*,\lambda^*)+\frac{\epsilon(t)}{2}\Vert x(t)\Vert^2)\\
			&+\frac{1}{2}\Vert c(t)(x(t)-x^*)+\dot{x}(t)\Vert ^2+\frac{b(t)}{2}\Vert \lambda(t)-\lambda^*\Vert^2+\frac{d(t)}{2}\Vert x(t)-x^*\Vert^2,
		\end{aligned}
	\end{equation}
	where $b(t)$ := $\frac{1}{\theta t^2}$, $c(t)$ := $\frac{1}{\theta t}$, $d(t)$ := $\frac{\theta t\gamma(t)-\theta-1}{\theta^2t^2}$.  We will provide specific conditions later to ensure that  $\widetilde{\mathcal{G}}(t)$ is a Lyapunov function for the dynamical system \eqref{equ5}.
	
	\begin{lemma}\label{lem2.2}
		Let $\mathcal{X},\mathcal{Y}$ be two real Hilbert spaces, $f:\mathcal{X}\rightarrow\mathbb{R}$ be a continuously differentiable convex function such that $\nabla f$ is Lipschitz continuous,  $A:\mathcal{X}\rightarrow\mathcal{Y}$ be a continuous linear operator, $b\in\mathcal{Y}$ be a given point in $\mathcal{Y}$, and $\theta, t_0$ be two positive constants. Assume that $\gamma,\beta,\epsilon:\left[t_0,+\infty\right)\rightarrow\left[0,+\infty\right)$ are three differentiable functions. Then, for any trajectory (global solution) $(x(t),\lambda(t))$ of the dynamical system $\rm{\eqref{equ5}}$ and any primal-dual optimal solution $(x^*,\lambda^*)\in\Omega$ of the linear equality constrained optimization problem $\rm{\eqref{equ1}}$,
		\begin{equation}\label{equ37}
			\begin{aligned}
				\frac{2}{t}\widetilde{\mathcal{G}}(t)+\dot{\widetilde{\mathcal{G}}}(t)\leq\;&(\dot{\beta}(t)+\frac{2\theta-1}{\theta t}\beta(t))(\mathcal{L}_{\sigma}(x(t),\lambda^*)-\mathcal{L}_{\sigma}(x^*,\lambda^*))\\
				&+\frac{1}{2}(\dot{\beta}(t)\epsilon(t)+\beta(t)\dot{\epsilon}(t)+\frac{2\theta-1}{\theta t}\beta(t)\epsilon(t))\Vert x(t)\Vert^2\\
				&+\frac{1}{2\theta t^2}(\gamma(t)+t\dot{\gamma}(t)-t\beta(t)\epsilon(t))\Vert x(t)-x^*\Vert^2\\
				&-\frac{\sigma}{2\theta t}\beta(t)\Vert Ax(t)-b\Vert^2+(\frac{1+\theta }{\theta t}-\gamma(t))\Vert\dot{x}(t)\Vert^2+\frac{1}{2\theta t}\beta(t)\epsilon(t)\Vert x^*\Vert^2.
			\end{aligned}
		\end{equation}
	\end{lemma}
	{\it Proof} 
	Calculating the derivative of \eqref{equ36} yields
	\begin{equation*}
		\begin{aligned}
			\dot{\widetilde{\mathcal{G}}}(t) =\;&\dot{\beta}(t)(\mathcal{L}_{\sigma}(x(t),\lambda^*)-\mathcal{L}_{\sigma}(x^*,\lambda^*))+\beta(t)\langle\nabla_x\mathcal{L}_{\sigma}(x(t),\lambda^*),\dot{x}(t)\rangle+c(t)\Vert\dot{x}(t)\Vert^2\\
			&+\frac{1}{2}(\dot{\beta}(t)\epsilon(t)+\beta(t)\dot{\epsilon}(t))\Vert x(t)\Vert^2+\beta(t)\epsilon(t)\langle x(t),\dot{x}(t)\rangle+(c(t)\dot{c}(t)+\frac{\dot{d}(t)}{2})\Vert x(t)-x^*\Vert^2\\
			&+(c(t)^2+\dot{c}(t)+d(t))\langle x(t)-x^*,\dot{x}(t)\rangle+\frac{\dot{b}(t)}{2}\Vert \lambda(t)-\lambda^*\Vert^2+b(t)\langle \lambda(t)-\lambda^*,\dot{\lambda}(t)\rangle\\
			&+\langle c(t)(x(t)-x^*)+\dot{x}(t),\ddot{x}(t)\rangle,
		\end{aligned}
	\end{equation*}
	which together with the dynamical system \eqref{equ5} indicates that
	\begin{equation*}
		\begin{aligned}
			\dot{\widetilde{\mathcal{G}}}(t) =\;&\dot{\beta}(t)(\mathcal{L}_{\sigma}(x(t),\lambda^*)-\mathcal{L}_{\sigma}(x^*,\lambda^*))+\frac{1}{2}(\dot{\beta}(t)\epsilon(t)+\beta(t)\dot{\epsilon}(t))\Vert x(t)\Vert^2+(c(t)\dot{c}(t)+\frac{\dot{d}(t)}{2})\Vert x(t)-x^*\Vert^2\\
			&+(c(t)^2+\dot{c}(t)+d(t)-c(t)\gamma(t))\langle x(t)-x^*,\dot{x}(t)\rangle+(c(t)-\gamma(t))\Vert\dot{x}(t)\Vert^2+\frac{\dot{b}(t)}{2}\Vert \lambda(t)-\lambda^*\Vert^2\\
			&+(\theta t^2b(t)\beta(t)-\beta(t))\langle \lambda(t)-\lambda^*, A\dot{x}(t)\rangle+(tb(t)\beta(t)-c(t)\beta(t))\langle\lambda(t)-\lambda^*, Ax(t)-b\rangle\\
			&-c(t)\beta(t)\epsilon(t)\langle x(t)-x^*,x(t)\rangle-c(t)\beta(t)\langle \nabla_x\mathcal{L}_{\sigma}(x(t),\lambda^*),x(t)-x^*\rangle.
		\end{aligned}
	\end{equation*}
	Recall the estimation \eqref{equ10} and thus we can get that 
	\begin{equation}\label{eq-4}
		\begin{aligned}
			\dot{\widetilde{\mathcal{G}}}(t) \leq\;&(\dot{\beta}(t)-c(t)\beta(t))(\mathcal{L}_{\sigma}(x(t),\lambda^*)-\mathcal{L}_{\sigma}(x^*,\lambda^*))+\frac{1}{2}(\dot{\beta}(t)\epsilon(t)+\beta(t)\dot{\epsilon}(t)-c(t)\beta(t)\epsilon(t))\Vert x(t)\Vert^2\\
			&+(c(t)\dot{c}(t)+\frac{\dot{d}(t)}{2}-\frac{1}{2}c(t)\beta(t)\epsilon(t))\Vert x(t)-x^*\Vert^2+(c(t)^2+\dot{c}(t)+d(t)-c(t)\gamma(t))\langle x(t)-x^*,\dot{x}(t)\rangle\\
			&+(c(t)-\gamma(t))\Vert\dot{x}(t)\Vert^2+\frac{1}{2}c(t)\beta(t)\epsilon(t)\Vert x^*\Vert^2+\frac{\dot{b}(t)}{2}\Vert\lambda(t)-\lambda^*\Vert^2-\frac{\sigma}{2}c(t)\beta(t)\Vert Ax(t)-b\Vert^2\\
			&+(\theta t^2b(t)\beta(t)-\beta(t))\langle \lambda(t)-\lambda^*, A\dot{x}(t)\rangle+(tb(t)\beta(t)-c(t)\beta(t))\langle\lambda(t)-\lambda^*, Ax(t)-b\rangle.
		\end{aligned}
	\end{equation}
	Then, it follows from \eqref{equ36}  and \eqref{eq-4} that 
	\begin{equation}\label{equ46}
		\begin{aligned}
			\frac{2}{t}\widetilde{\mathcal{G}}(t)+\dot{\widetilde{\mathcal{G}}}(t)\leq\;&(\dot{\beta}(t)-c(t)\beta(t)+\frac{2}{t}\beta(t))(\mathcal{L}_{\sigma}(x(t),\lambda^*)-\mathcal{L}_{\sigma}(x^*,\lambda^*))\\
			&+\frac{1}{2}(\dot{\beta}(t)\epsilon(t)+\beta(t)\dot{\epsilon}(t)-(c(t)-\frac{2}{t})\beta(t)\epsilon(t))\Vert x(t)\Vert^2-\frac{\sigma}{2}c(t)\beta(t)\Vert Ax(t)-b\Vert^2\\
			&+(c(t)-\gamma(t)+\frac{1}{t})\Vert\dot{x}(t)\Vert^2+\frac{1}{2}c(t)\beta(t)\epsilon(t)\Vert x^*\Vert^2+(\frac{\dot{b}(t)}{2}+\frac{b(t)}{t})\Vert\lambda(t)-\lambda^*\Vert^2\\
			&+(c(t)\dot{c}(t)+\frac{\dot{d}(t)}{2}-\frac{1}{2}c(t)\beta(t)\epsilon(t)+\frac{1}{t}c(t)^2+\frac{d(t)}{t})\Vert x(t)-x^*\Vert^2\\
			&+(\theta t^2b(t)\beta(t)-\beta(t))\langle \lambda(t)-\lambda^*, A\dot{x}(t)\rangle\\
			&+(tb(t)\beta(t)-c(t)\beta(t))\langle\lambda(t)-\lambda^*, Ax(t)-b\rangle\\
			&+ (c(t)^2+\dot{c}(t)+d(t)-c(t)\gamma(t)+\frac{2}{t}c(t))\langle x(t)-x^*,\dot{x}(t)\rangle.
		\end{aligned}
	\end{equation}
	According to the specified formulations of $b(\cdot)$, $c(\cdot)$ and $d(\cdot)$, we can obtain that 
	\begin{equation*}
		\begin{cases}
			c(t)^2+\dot{c}(t)+d(t)-c(t)\gamma(t)+\frac{2}{t}c(t)=0,\\
			\theta t^2b(t)\beta(t)-\beta(t)=0,\\
			tb(t)\beta(t)-c(t)\beta(t)=0,\\
			\frac{\dot{b}(t)}{2}+\frac{b(t)}{t}=0.
		\end{cases}
	\end{equation*}
	Therefore, \eqref{equ46} implies that \eqref{equ37} holds. This completes the proof of Lemma \ref{lem2.2}.
	\qed
	
	\begin{theorem}\label{thm2.2}
		Let all hypotheses in Lemma \ref{lem2.2} hold. Assume further that $\epsilon$ is $\mathcal{C}^1$ and non-increasing function such that $\int_{t_0}^{+\infty}\frac{\beta(t)\epsilon(t)}{t} dt<+\infty$,  and $\beta, \gamma$ are two $\mathcal{C}^1$ and  non-negative functions satisfying  $\rm{\eqref{equ6}}$-${\rm\eqref{equ8}}$. Then, for any trajectory $(x(t),\lambda(t))$ of the dynamical system $\rm{\eqref{equ5}}$ and any $(x^*,\lambda^*)\in\Omega$, the following conclusions hold:
		\begin{flalign*}
			&(a)\;\lim\limits_{t\rightarrow+\infty}\Vert \frac{1}{\theta t}(x(t)-x^*)+\dot{x}(t)\Vert=0;&\\
			&(b)\;\text{ When }\liminf\limits_{t\rightarrow +\infty}\beta(t)\neq 0,\;\lim\limits_{t\rightarrow+\infty}\mathcal{L}_{\sigma}(x(t),\lambda)-\mathcal{L}_{\sigma}(x^*,\lambda^*)=0;&\\
			&\quad\;\;\text{Specially, when }\lim\limits_{t\rightarrow+\infty}\beta(t)=+\infty,&\\
			&\qquad\quad\mathcal{L}_{\sigma}(x(t),\lambda)-\mathcal{L}_{\sigma}(x^*,\lambda^*)=o(\frac{1}{\beta(t)})\quad\text{and}\quad\Vert\nabla f(x(t))-\nabla f(x^*)\Vert = o(\frac{1}{\sqrt{\beta(t)}}).&
		\end{flalign*}
	\end{theorem}
	{\it Proof} 
	According to the assumptions on  the coefficients $\epsilon(t)$, $\beta(t)$, $\gamma(t)$ in Theorem \ref{thm2.2}, it follows  that
	\begin{equation*}
		\begin{cases}
			\vspace{2mm}
			\dot{\beta}(t)+\frac{2\theta-1}{\theta t}\beta(t)\leq  0,\\
			
			\vspace{2mm}
			\frac{1}{2\theta t^2}(\gamma(t)+t\dot{\gamma}(t)-t\beta(t)\epsilon(t))\leq 0,\\
			\vspace{2mm}
			\frac{1+\theta }{\theta t}-\gamma(t)\leq 0,\\
			\vspace{2mm}
			\dot{\beta}(t)\epsilon(t)+\beta(t)\dot{\epsilon}(t)+\frac{2\theta-1}{\theta t}\beta(t)\epsilon(t)\leq 0,\\
			\vspace{2mm}
			-\frac{\sigma}{2\theta t}\beta(t)\leq 0,
		\end{cases}
	\end{equation*}
	where the first to third inequalities hold directly from \eqref{equ6}-\eqref{equ8} and the non-negativity of $\theta$, the fourth inequality is satisfied since $\epsilon(t)$ is a non-negative and non-increasing function, and $\beta(t)$ is a non-negative function satisfying \eqref{equ6}, and the last inequality holds due to the non-negativity of $\sigma$, $\theta$ and $\beta(t)$. This together with estimation \eqref{equ37} implies that for all $t\ge t_0$,
	\begin{equation}\label{equ21}
		\frac{2}{t}\widetilde{\mathcal{G}}(t)+\dot{\widetilde{\mathcal{G}}}(t)\leq\frac{1}{2\theta t}\beta(t)\epsilon(t)\Vert x^*\Vert^2.
	\end{equation}
	Multiplying both sides of \eqref{equ21} by $t^2$ and integrating the obtained results on $[t_0,t]$ yield
	\begin{equation}\label{equ23}
		\widetilde{\mathcal{G}}(t)\leq \frac{t_0^2}{t^2}\widetilde{\mathcal{G}}(t_0)+\frac{\Vert x^*\Vert^2}{2\theta t^2}\int_{t_0}^{t}s\beta(s)\epsilon(s)ds.
	\end{equation}
	This allows us to use Lemma\ref{lemB.2} with $\psi(t)=t^2$ and $\phi(t)=\frac{\beta(t)\epsilon(t)}{t}$ to \eqref{equ23} to get that 
	\begin{equation}\label{eq-5}
		\lim\limits_{t\rightarrow+\infty}\frac{1}{t^2}\int_{t_0}^{t}s\beta(s)\epsilon(s)ds=0,
	\end{equation}
	where the condition $\int_{t_0}^{+\infty}\frac{\beta(t)\epsilon(t)}{t} dt<+\infty$ is used.  Since $\widetilde{\mathcal{G}}(t)\ge 0$, it follows from $\eqref{equ23}$ and $\eqref{eq-5}$ that $\lim\limits_{t\rightarrow+\infty}\widetilde{\mathcal{G}}(t)=0$ which indicates that $\widetilde{\mathcal{G}}(t)$ is a Lyapunov function and moreover the definition of $\widetilde{\mathcal{G}}(t)$ implies that
	\begin{equation*}
		\lim\limits_{t\rightarrow+\infty}\Vert \frac{1}{\theta t}(x(t)-x^*)+\dot{x}(t)\Vert=0 ,
	\end{equation*}
	and 
	\begin{equation*}
		\lim\limits_{t\rightarrow+\infty}\beta(t)(\mathcal{L}_{\sigma}(x(t),\lambda^*)-\mathcal{L}_{\sigma}(x^*,\lambda^*))= 0,
	\end{equation*}
	which implies
	\begin{equation*}
		\lim\limits_{t\rightarrow+\infty}\mathcal{L}_{\sigma}(x(t),\lambda^*)-\mathcal{L}_{\sigma}(x^*,\lambda^*)= 0,
	\end{equation*}
	when $\liminf\limits_{t\rightarrow+\infty}\beta(t)\neq 0$, and
	\begin{equation}\label{equ39}
		\mathcal{L}_{\sigma}(x(t),\lambda^*)-\mathcal{L}_{\sigma}(x^*,\lambda^*)=o(\frac{1}{\beta(t)}),
	\end{equation}
	when $\lim\limits_{t\rightarrow+\infty}\beta(t)=+\infty$.
	Recall the estimation  \eqref{equ38} and multiply at its both sides the non-negative $\beta(t)$ to get that
	\begin{equation*}
		\frac{\beta(t)}{2L}\Vert \nabla f(x(t))-\nabla f(x^*)\Vert^2\leq\beta(t)(\mathcal{L}_{\sigma}(x(t),\lambda^*)-\mathcal{L}_{\sigma}(x^*,\lambda^*)).
	\end{equation*}
	According to \eqref{equ39}, this implies
	\begin{equation*}
		\Vert\nabla f(x(t))-\nabla f(x^*)\Vert = o(\frac{1}{\sqrt{\beta(t)}}).
	\end{equation*}
	This completes the proof of Theorem \ref{thm2.2}.
	\qed
	
	\begin{remark}
		When the  linear equality constrained convex optimization problem \eqref{equ1} reduces to an unconstrained optimization problem and the parameter $\gamma(t)$ is set as $\frac{\alpha}{t}$, the dynamical system \eqref{equ5} degenerates to the inertial second-order dynamical system with Tikhonov regularization proposed by Xu et al. \cite{bib51}. In \cite{bib51}, Xu et al. studied separately the asymptotic convergence properties and strong convergence of the trajectory under the two distinct conditions:${\int_{t_0}^{+\infty}t\beta(t)\epsilon(t)  dt<+\infty}$ and ${\int_{t_0}^{+\infty}\frac{\beta(t)\epsilon(t)}{t}  dt=+\infty}$. In contrast, Theorem \ref{thm2.2} establishes the asymptotic convergence properties of the dynamical system under the condition ${\int_{t_0}^{+\infty}\frac{\beta(t)\epsilon(t)}{t}  dt<+\infty}$, thus filling the gap in the analysis of the inertial second-order dynamical system with Tikhonov regularization proposed in \cite{bib51} under this condition.
	\end{remark}
	
	\subsection{\bf Particular Cases}
	
	In this subsection, we consider some particular cases by choosing specific parameters $\theta$, $\beta(t)$ and $\gamma(t)$.  Specifically, we let $\beta(t)=t^\beta\;(\beta\ge 0)$ and consider the following two cases according the condition \eqref{equ7}.
	
	\vspace{0.5em}
	{\large \noindent {\bf Case 1:} $\bm{\gamma(t)+t\dot{\gamma}(t)\leq 0}$}
	\vspace{0.5em}
	
	\noindent  In this case, we consider the choice $\gamma(t)=\frac{\alpha}{t}$, which is commonly used in the literature for the dynamical systems of unconstrained/constrained optimization problems. Moreover, under the assumptions \eqref{equ6} and \eqref{equ8}, it holds that $\frac{1}{\alpha-1}\leq\theta\leq\frac{1}{\beta+2}$. Based on Theorems \ref{thm2.1} and \ref{thm2.2}, we then obtain the following results.
	\begin{theorem}\label{thm2.3}
		In the dynamical system $\rm{\eqref{equ5}}$, let
		\vspace{-7pt}
		\begin{equation*}
			\gamma(t)=\frac{\alpha}{t},\quad\beta(t)=t^{\beta},\quad \frac{1}{\alpha-1}\leq\theta\leq\frac{1}{\beta+2},
		\end{equation*}
		where $\alpha, \beta$ are constants with $\alpha\ge\beta+3$ and $ \beta\ge0$. Suppose that $t_0>0$, $\epsilon:\left[t_0,+\infty\right)\rightarrow\left[0,+\infty\right)$ is a $\mathcal{C}^1$ and non-increasing function, and  $(x(t),\lambda(t))$ with $t\geq t_0$ is a global solution of the dynamical system $\rm{\eqref{equ5}}$. 
		
		\noindent{\rm{(i)}}\;\; If $\int_{t_0}^{+\infty}t^{\beta+1}\epsilon(t) dt<+\infty$, then for any $(x^*,\lambda^*)\in\Omega$,
		\vspace{-5pt}
		\begin{adjustwidth}{5mm}{}
		\begin{flalign*}
		 &\begin{array}{l@{\quad}l}
			\textbf{(Boundedness of the Trajectory)} & \text{the trajectory }(x(\cdot),\lambda(\cdot))\text{ is bounded};\\[5pt]
			\textbf{(Pointwise Estimates)} & 
			\left\{
			\begin{array}{l}
				\Vert \dot{x}(t)\Vert=\mathcal{O}(\frac{1}{t}),\\[5pt]
				\Vert Ax(t)-b\Vert=\mathcal{O}(t^{-(\beta+2)}),\\[5pt]
				\vert f(x(t))-f(x^*)\vert =\mathcal{O}(t^{-(\beta+2)}),\\[5pt]
				\Vert  \nabla f(x(t))-\nabla f(x^*)\Vert=\mathcal{O}(t^{-\frac{\beta+2}{2}}),\\[5pt]
				\mathcal{L}_{\sigma}(x(t),\lambda^*)-\mathcal{L}_{\sigma}(x^*,\lambda^*)=\mathcal{O}(t^{-(\beta+2)});
			\end{array}
			\right. \\[40pt]
			\textbf{(Integral Estimates)} & 
			\begin{array}{l}	
				\int_{t_0}^{+\infty}t^{\beta+1}\Vert Ax(t)-b\Vert^2 dt<+\infty;
			\end{array}
		\end{array}&
		\end{flalign*}
			\end{adjustwidth}
		
			\qquad In particular, if $\frac{1}{\alpha-1}\leq\theta<\frac{1}{\beta+2}$, then 
			\begin{equation*}
				\begin{cases}
					\int_{t_0}^{+\infty} t^{\beta+1}(\mathcal{L}_{\sigma}(x(t),\lambda^*)-\mathcal{L}_{\sigma}(x^*,\lambda^*))dt<+\infty,\\[5pt]
					\int_{t_0}^{+\infty} t^{\beta+1} \Vert  \nabla f(x(t))-\nabla f(x^*)\Vert^2 dt<+\infty,
				\end{cases}
			\end{equation*}
			
			\qquad and if $\frac{1}{\alpha-1}<\theta\leq\frac{1}{\beta+2}$, then $	\int_{t_0}^{+\infty}t\Vert \dot{x}(t)\Vert dt<+\infty$;
		\\
		{\rm{(ii)}}\;\; If  $\int_{t_0}^{+\infty}t^{\beta-1}\epsilon(t)dt<+\infty$, then for any $(x^*,\lambda^*)\in\Omega$, 
		\vspace{-10pt}
		\begin{adjustwidth}{5mm}{}
			\begin{flalign*}
				&(a)\;\lim\limits_{t\rightarrow+\infty}\mathcal{L}_{\sigma}(x(t),\lambda)-\mathcal{L}_{\sigma}(x^*,\lambda^*)=0; \text{ Specially, when }\beta>0,&\\
				&\quad\;\;\mathcal{L}_{\sigma}(x(t),\lambda)-\mathcal{L}_{\sigma}(x^*,\lambda^*)=o(t^{-\beta})\vspace{1mm}\quad \text{and}\quad\Vert\nabla f(x(t))-\nabla f(x^*)\Vert = o(t^{-\frac{\beta}{2}});&\\
				&(b)\;\lim\limits_{t\rightarrow+\infty}\Vert \frac{1}{\theta t}(x(t)-x^*)+\dot{x}(t)\Vert=0.&
			\end{flalign*}
		\end{adjustwidth}

	\end{theorem}
	
	{\it Proof} 
	By simple calculation, it is easy to check that 
	$$(2\theta-1)t\beta(t)+\theta t^2\dot{\beta}(t)=((2+\beta)\theta-1)t^{\beta+1},$$
	and 
	$$\gamma(t)+t\dot{\gamma}(t)=0,$$
	which together with $\frac{1}{\alpha-1}\leq\theta\leq\frac{1}{\beta+2}$ implies the parameters  satisfy the conditions \eqref{equ6}-\eqref{equ8} in Theorem \ref{thm2.1} and all conditions in Theorem \ref{thm2.2} when $t\ge t_0>0$. In particular, if $\frac{1}{\alpha-1}\leq\theta<\frac{1}{\beta+2}$, then  in combination with $\beta(t)=t^\beta$, it is straightforward to verify that $(2\theta-1)\beta(t)+\theta t\dot{\beta}(t)< 0,\;\forall t\ge t_0$. Moreover, if $\frac{1}{\alpha-1}<\theta\leq\frac{1}{\beta+2}$, then together with $\gamma(t)=\frac{\alpha}{t}$, one can similarly check that $\theta t\gamma(t)-\theta-1>0,\;\forall t\ge t_0$. Therefore, Theorem \ref{thm2.3} is a direct consequence of  Theorem \ref{thm2.1} and Theorem \ref{thm2.2}. This completes the proof of Theorem \ref{thm2.3}.
	\qed
	
	\begin{remark}
		(1) In particular, if $\beta = 0 \;(i.e., \beta(t)=1)$, $\gamma(t) = \frac{\alpha}{t}$ and $\theta = \frac{1}{\alpha-1}$,  Theorem \ref{thm2.3} reduces to  Theorem 4.1 and Theorem 4.6 in {\rm\cite{bib23}}. (2) There are other choices for this case, for example $\gamma(t)=\frac{1+\alpha t}{t^2}$ with $\alpha>1$. In this case, we can further obtain $\int_{t_0}^{+\infty}((\alpha\theta-\theta-1)t+\theta)\Vert \dot{x}(t)\Vert dt<+\infty$ if $\theta\in\left[\frac{1}{\alpha-1}, \frac{1}{\beta+2}\right]$.
	\end{remark}

	\noindent {\large {\bf Case 2:} $\bm{0\leq\gamma(t)+t\dot{\gamma}(t)\leq t\beta(t)\epsilon(t)}$}
	
	\vspace{ 2mm}
	\noindent  In this case, we let $\gamma(t)=\frac{2\alpha t-1}{t^2}$,  $\epsilon(t)=\frac{1}{t^{\beta+3}}$,   and, by the assumptions \eqref{equ6} and \eqref{equ8}, we have $\frac{1}{2\alpha-2}\leq\theta\leq\frac{1}{\beta+2}$. Thus we can get from Theorems \ref{thm2.1} and \ref{thm2.2} the following theorem. 
	
	\begin{theorem}\label{thm2.5}		
		Suppose that the parameters of dynamical system $\rm{\eqref{equ5}}$ satisfy
		\begin{equation}\label{eq-1-7}
			\gamma(t)=\frac{2\alpha t-1}{t^2},\;\beta(t)=t^\beta,\;\epsilon(t)=\frac{1}{t^{\beta+3}},\;\frac{1}{2\alpha-2}\leq\theta\leq\frac{1}{\beta+2},
		\end{equation} 
		where $\alpha, \beta$ are constants with $\alpha>\frac{\beta+4}{2}$ and $\beta\ge0$. Suppose that $t_0\ge 1$, $\epsilon:\left[t_0,+\infty\right)\rightarrow\mathbb{R}_+$ is a $\mathcal{C}^1$ and non-increasing function and $(x(t),\lambda(t))$ with $t\geq t_0$ is a global solution of dynamical system $\rm{\eqref{equ5}}$. Then, for any $(x^*,\lambda^*)\in\Omega$, the following conclusions hold:
		\vspace{-7pt}
		\[
		\begin{array}{l}
			\begin{array}{ll}
				(1) \textbf{ (Boundedness of Trajectory)} & \text{the trajectory }(x(\cdot),\lambda(\cdot))\text{ is bounded};\\[8pt]
				(2) \textbf{ (Pointwise Estimates)} & 
				\left\{
				\begin{array}{l}
					\Vert \dot{x}(t)\Vert=\mathcal{O}(\frac{1}{t}),\\[5pt]
					\Vert Ax(t)-b\Vert=\mathcal{O}(t^{-(\beta+2)}),\\[5pt]
					\vert f(x(t))-f(x^*)\vert =\mathcal{O}(t^{-(\beta+2)}),\\[5pt]
					\Vert  \nabla f(x(t))-\nabla f(x^*)\Vert=\mathcal{O}(t^{-\frac{\beta+2}{2}}),\\[5pt]
					\mathcal{L}_{\sigma}(x(t),\lambda^*)-\mathcal{L}_{\sigma}(x^*,\lambda^*)=\mathcal{O}(t^{-(\beta+2)}),\\[5pt]
					\lim\limits_{t\rightarrow+\infty}\Vert \frac{1}{\theta t}(x(t)-x^*)+\dot{x}(t)\Vert=0;
				\end{array}
				\right. \\[50pt]
				(3) \textbf{ (Integral Estimates)} & 
				\begin{array}{l}	
					\int_{t_0}^{+\infty}t^{\beta+1}\Vert Ax(t)-b\Vert^2 dt<+\infty;
				\end{array}
			\end{array}
			\vspace{2mm}
			\\[50pt]
			
			\begin{array}{l}	
					\text{\qquad In particular, if }\frac{1}{2\alpha-2}\leq\theta<\frac{1}{\beta+2}, \text{ then } 
						\begin{cases}
							\int_{t_0}^{+\infty} t^{\beta+1}(\mathcal{L}_{\sigma}(x(t),\lambda^*)-\mathcal{L}_{\sigma}(x^*,\lambda^*))dt<+\infty,\\[5pt]
							\int_{t_0}^{+\infty} t^{\beta+1} \Vert  \nabla f(x(t))-\nabla f(x^*)\Vert^2 dt<+\infty,
						\end{cases}
				\\[15pt]
				\text{\qquad and if } \frac{1}{2\alpha-2}<\theta\leq\frac{1}{\beta+2}, \text{ then } \int_{t_0}^{+\infty}((2\alpha\theta-\theta-1)t-\theta )\Vert \dot{x}(t)\Vert dt<+\infty;
					\\[10pt]\vspace{2mm}
				(4) \;\lim\limits_{t\rightarrow+\infty}\mathcal{L}_{\sigma}(x(t),\lambda)-\mathcal{L}_{\sigma}(x^*,\lambda^*)=0; \text{ Specially, when }\beta>0,\\[5pt]
				\quad\;\;\;\mathcal{L}_{\sigma}(x(t),\lambda)-\mathcal{L}_{\sigma}(x^*,\lambda^*)=o(t^{-\beta})\vspace{1mm}\quad\text{and}\quad\Vert\nabla f(x(t))-\nabla f(x^*)\Vert = o(t^{-\frac{\beta}{2}});\\[5pt]
			\end{array}
			
		\end{array}
		\]
		
	\end{theorem}
	{\it Proof} With specific formulations of parameters in \eqref{eq-1-7}, we can obtain by simple calculating that 
	\begin{itemize}
		\item[\quad ] \hspace{4mm} $\gamma(t)+t\dot{\gamma}(t)=\frac{1}{t^2}=t\beta(t)\epsilon(t)$,\\
		\item[\quad ] \hspace{4mm}  $\int_{t_0}^{+\infty}t\beta(t)\epsilon(t)dt=\frac{1}{t_0}$ and  $\int_{t_0}^{+\infty}\frac{\beta(t)\epsilon(t)}{t}dt=\frac{1}{3t_0^3}$,\\
		\item[\quad ] \hspace{4mm} $(2\theta-1)t\beta(t)+\theta t^2\dot{\beta}(t)=((2+\beta)\theta-1)t^{\beta+1}\leq 0$,\\
		\item[\quad ] \hspace{4mm}  $\theta t \gamma(t)-\theta-1\ge 0$,
	\end{itemize}
	which indicates that the conditions \eqref{equ6}-\eqref{equ8} in Theorem \ref{thm2.1} and all conditions in Theorem \ref{thm2.2} hold for all $t\ge t_0>1$. In particular, if $\frac{1}{2\alpha-2}\leq\theta<\frac{1}{\beta+2}$, then  in combination with $\beta(t)=t^\beta$, it is straightforward to verify that $(2\theta-1)\beta(t)+\theta t\dot{\beta}(t)< 0,\;\forall t\ge t_0>1$. Moreover, if $\frac{1}{2\alpha-2}<\theta\leq\frac{1}{\beta+2}$, then together with $	\gamma(t)=\frac{2\alpha t-1}{t^2}$, one can similarly check that $\theta t\gamma(t)-\theta-1>0,\;\forall t\ge t_0>1$. Therefore, according to Theorem \ref{thm2.1} and Theorem \ref{thm2.2}, we can obtain the conclusions in Theorem \ref{thm2.5}. This completes the proof of Theorem \ref{thm2.5}.
	\qed
	%

	
	%
	\section{Strong Convergence of the Trajectory }\label{sec3}
	In this section,   when the Tikhonov regularization parameter $\epsilon(t)$ approaches zero at an appropriate rate, we provide a convergence result for the trajectory of the dynamical system \eqref{equ5}, i.e., the trajectory generated by the dynamical system \eqref{equ5} strongly converges to the minimum norm solution of the problem \eqref{equ1}.
	
	\subsection{\bf Classical Facts}\label{subsec3}
	We denote the unique element of minimal norm of the solution set $S$ of the problem \eqref{equ1} by $\hat{x}^*$. Then there exists an optimal dual solution $\hat{\lambda}^*\in\mathcal{Y}$ of the Lagrange dual problem $\max_{\lambda}\min_{x\in\mathcal{X}}\widetilde{\mathcal{L}}(x,\lambda)$ such that $(\hat{x}^*,\hat{\lambda}^*)\in \Omega$, where $\widetilde{\mathcal{L}}:\mathcal{X}\times\mathcal{Y}\rightarrow\mathbb{R}$ is a Lagrangian function associated with the problem \eqref{equ1}.
	And according to the definition of $\hat{x}^*$, it has a remarkable geometrical property
	\begin{equation*}
		\hat{x}^*= \rm{Proj}_S 0.
	\end{equation*}
	For any $\epsilon>0$, we denote by $x_\epsilon$ the unique solution of the strongly convex minimization problem
	\begin{equation*}
		x_\epsilon:={\rm{arg}}\min_{x\in \mathcal{X}}\left\{\mathcal{L}_{\sigma}(x,\hat{\lambda}^*)+\frac{\epsilon}{2}\Vert x\Vert^2\right\}.
	\end{equation*}
	The first-order optimality condition leads to
	\begin{equation}\label{equ41}
		\nabla_x\mathcal{L}_{\sigma}(x_\epsilon,\hat{\lambda}^*)+\epsilon x_\epsilon=0.
	\end{equation}
	Owing to the classical properties of the Tikhonov regularization, it follows that
	\begin{equation*}
		\lim\limits_{\epsilon\rightarrow 0}\Vert x_\epsilon-\hat{x}^*\Vert =0.
	\end{equation*}
	The result was first obtained by Tikhonov \cite{bib31} for optimization problems and Browder \cite{bib33} for monotone variational inequalities. We refer to \cite{bib19} and \cite{bib25} for more references. Moreover, based on the convexity of $\mathcal{L}_{\sigma}(x, \hat{\lambda}^*)$, we know that $\nabla_x\mathcal{L}_{\sigma}(x,\hat{\lambda}^*)$ is a monotone operator. By using $\nabla_x\mathcal{L}_{\sigma}(\hat{x}^*, \hat{\lambda}^*)=0$ and \eqref{equ41},  one arrives at
	\begin{equation*}
		\langle x_\epsilon-\hat{x}^*,-\epsilon x_\epsilon\rangle\ge 0,
	\end{equation*}
	which indicates that
	\begin{equation*}
		\Vert x_\epsilon\Vert\leq\Vert \hat{x}^*\Vert,\; \forall\epsilon>0.
	\end{equation*}
	
	\subsection{\bf Strong Convergence }\label{subsec4}
	\begin{theorem}\label{thm3.1}
		Supposed that $f:\mathcal{X}\rightarrow\mathbb{R}$ is a continuously differentiable convex function such that $\nabla f$ is Lipschitz continuous. Let $\theta>0$, $t_0>0$ and $\epsilon:\left[t_0,+\infty\right)\rightarrow\left[0,+\infty\right)$ is a $\mathcal{C}^1$ and non-increasing function, and let $\gamma,\beta:\left[t_0,+\infty\right)\rightarrow\left[0,+\infty\right)$ be two $\mathcal{C}^1$ and  non-negative functions satisfying $\rm{\eqref{equ6}}$-$\rm{\eqref{equ8}}$, $\int_{t_0}^{+\infty}\frac{\beta(t)\epsilon(t)}{t} dt<+\infty$, $\liminf\limits_{t\rightarrow +\infty}\beta(t)\neq 0$ and $\lim\limits_{t\rightarrow+\infty}t^2\beta(t)\epsilon(t)=+\infty$.
		For the global solution $(x(\cdot),\lambda(\cdot))$ of the dynamical system $\rm{\eqref{equ5}}$, we have 
		\begin{equation*}
			\liminf_{t\rightarrow+\infty}\Vert x(t)-\hat{x}^*\Vert=0,
		\end{equation*}
		where $\hat{x}^*=\rm{Proj}_S 0$ is the unique element of minimal norm in $S$. Moreover, if there exists a large enough T such that the trajectory $\left\{x(t):t\ge T\right\}$ stays in either the open ball $B(0,\Vert \hat{x}^*\Vert)$ or its complement, then 
		\begin{equation*}\label{equ24}
			\lim\limits_{t\rightarrow+\infty}\Vert x(t)-\hat{x}^*\Vert=0.
		\end{equation*}
	\end{theorem}
	{\it Proof} We prove the theorem by considering the following three cases relative to the trajectory $x(\cdot)$.
	
	\textbf{Case I:} There exists a large enough T such that $x(t) $ stays in the complement of $B(0,\Vert \hat{x}^*\Vert)$, i.e.,
	\begin{equation}\label{equ25}
		\Vert x(t)\Vert\ge\Vert \hat{x}^*\Vert,\; \forall t\ge T.
	\end{equation}
	Define a  function $\widehat{\mathcal{G}}:\left[t_0,+\infty\right)\rightarrow\left[0,+\infty\right)$ as follows 
	\begin{equation*}\label{equ26}
		\begin{aligned}
			\widehat{\mathcal{G}}(t) :=\;&\beta(t)(\mathcal{L}_{\sigma}(x(t),\hat{\lambda}^*)-\mathcal{L}_{\sigma}(\hat{x}^*,\hat{\lambda}^*)+\frac{\epsilon(t)}{2}\Vert x(t)\Vert^2-\frac{\epsilon(t)}{2}\Vert\hat{x}^*\Vert^2)+\frac{1}{2}\Vert \frac{1}{\theta t}(x(t)-\hat{x}^*)+\dot{x}(t)\Vert ^2\\
			&+\frac{\theta t\gamma(t)-\theta-1}{2\theta^2t^2}\Vert x(t)-\hat{x}^*\Vert^2+\frac{1}{2\theta t^2}\Vert \lambda(t)-\hat{\lambda}^*\Vert^2,
		\end{aligned}
	\end{equation*}
	and calculate its derivative to get that
	\begin{equation*}
		\begin{aligned}
			\dot{\widehat{\mathcal{G}}}(t) =\;& \dot{\beta}(t)(\mathcal{L}_{\sigma}(x(t),\hat{\lambda}^*)-\mathcal{L}_{\sigma}(\hat{x}^*,\hat{\lambda}^*)+\frac{\epsilon(t)}{2}\Vert x(t)\Vert^2-\frac{\epsilon(t)}{2}\Vert\hat{x}^*\Vert^2)\\
			&+\frac{\beta(t)\dot{\epsilon}(t)}{2}(\Vert x(t)\Vert^2-\Vert \hat{x}^*\Vert^2)+\beta(t)\langle\nabla_x\mathcal{L}_{\sigma}(x(t),\hat{\lambda}^*),\dot{x}(t)\rangle+\beta(t)\epsilon(t)\langle x(t),\dot{x}(t)\rangle\\
			&+\frac{1}{2\theta t^3}(-t\gamma(t)+t^2\dot{\gamma}(t)+2)\Vert x(t)-\hat{x}^*\Vert^2+\frac{1}{\theta t}\Vert\dot{x}(t)\Vert^2-\frac{1}{\theta t^3}\Vert \lambda(t)-\hat{\lambda}^*\Vert ^2\\
			&+\langle\frac{1}{\theta t}(x(t)-\hat{x}^*)+\dot{x}(t),\ddot{x}(t)\rangle+\frac{1}{\theta t^2}(t\gamma(t)-2)\langle x(t)-\hat{x}^*,\dot{x}(t)\rangle+\frac{1}{\theta t^2}\langle \lambda(t)-\hat{\lambda}^*,\dot{\lambda}(t)\rangle.
		\end{aligned}
	\end{equation*}
	With similar arguments in the proof of Theorem \ref{thm2.2}, we can further obtain that
	\begin{equation*}
		\begin{aligned}
			\frac{2}{t}\widehat{\mathcal{G}}(t)+\dot{\widehat{\mathcal{G}}}(t)\leq&(\dot{\beta}(t)+\frac{2\theta -1}{\theta t}\beta(t))(\mathcal{L}_{\sigma}(x(t),\hat{\lambda}^*)-\mathcal{L}_{\sigma}(\hat{x}^*,\hat{\lambda}^*))\\
			&+\frac{1}{2\theta t^2}(\gamma(t)+t\dot{\gamma}(t)-t\beta(t)\epsilon(t))\Vert x(t)-\hat{x}^*\Vert^2+(\frac{1+\theta}{\theta t}-\gamma(t))\Vert\dot{x}(t)\Vert^2\\
			&-\frac{\sigma\beta(t)}{2\theta t}\Vert Ax(t)-b\Vert^2+\frac{1}{2}(\dot{\beta}(t)\epsilon(t)+\beta(t)\dot{\epsilon}(t)+\frac{2\theta-1}{\theta t}\beta(t)\epsilon(t))(\Vert x(t)\Vert^2-\Vert \hat{x}^*\Vert^2).
		\end{aligned}
	\end{equation*}
	Based on the assumptions on the parameters $\epsilon(t)$, $\beta(t)$, $\gamma(t)$,  the conditions \eqref{equ6}-\eqref{equ8}, and inequality  \eqref{equ25}, it follows that
	\begin{equation*}
		\frac{2}{t}\widehat{\mathcal{G}}(t)+\dot{\widehat{\mathcal{G}}}(t)\leq 0,\;\forall t\ge T.
	\end{equation*}
	Multiplying both sides of the above inequality by $t^2$  and integrating the obtained inequality  form $T$ to $t$ yield
	\begin{equation*}
		\widehat{\mathcal{G}}(t)\leq\frac{T^2\widehat{\mathcal{G}}(T)}{t^2},\;\forall t\ge T,
	\end{equation*}
	which indicates that 
	\begin{equation}\label{equ40}
		\beta(t)(\mathcal{L}_{\epsilon(t)}(x(t)-\mathcal{L}_{\epsilon(t)}(\hat{x}^*))\leq \widehat{\mathcal{G}}(t)\leq \frac{T^2\widehat{\mathcal{G}}(T)}{t^2},\;\forall t\ge T.
	\end{equation}
	Let $\hat{\lambda}^*$ be such that $(\hat{x}^*,\hat{\lambda}^*)\in\Omega$ and for each $\epsilon>0$, we define the function $\mathcal{L}_\epsilon:\mathcal{X}\rightarrow\mathbb{R}$ as
	\begin{equation*}
		\mathcal{L}_\epsilon(x) := \mathcal{L}_{\sigma}(x,\hat{\lambda}^*)+\frac{\epsilon}{2}\Vert x\Vert^2,
	\end{equation*}
	which is obviously strongly convex and  together with \eqref{equ41} results in
	\begin{equation}\label{equ27}
		\nabla\mathcal{L}_{\epsilon}(x_{\epsilon})=\nabla_x\mathcal{L}_{\sigma}(x_{\epsilon},\hat{\lambda}^*)+\epsilon x_{\epsilon}=0.
	\end{equation}
	Thus, for any $t\geq T$, we can get by the strong convexity of $\mathcal{L}_{\epsilon(t)}$ that
	\begin{equation*}
		\mathcal{L}_{\epsilon(t)}(x(t))-\mathcal{L}_{\epsilon(t)}(x_{\epsilon(t)})\ge \langle\nabla\mathcal{L}_{\epsilon(t)}(x_{\epsilon(t)}),x(t)-x_{\epsilon(t)}\rangle+\frac{\epsilon(t)}{2}\Vert x(t)-x_{\epsilon(t)}\Vert^2,
	\end{equation*}
	which implies by \eqref{equ27} that
	\begin{equation*}
		\mathcal{L}_{\epsilon(t)}(x(t))-\mathcal{L}_{\epsilon(t)}(x_{\epsilon(t)})\ge\frac{\epsilon(t)}{2}\Vert x(t)-x_{\epsilon(t)}\Vert^2.
	\end{equation*}
	Using the definition of $\mathcal{L}_{\epsilon}(x)$ and the fact that $ \mathcal{L}_{\sigma}(\hat{x}^*,\hat{\lambda}^*)\leq\mathcal{L}_{\sigma}(x_{\epsilon(t)},\hat{\lambda}^*)$, we derive
	\begin{equation*}
		\begin{aligned}
			\mathcal{L}_{\epsilon(t)}(\hat{x}^*)-\mathcal{L}_{\epsilon(t)}(x_{\epsilon(t)})=&\mathcal{L}_{\sigma}(\hat{x}^*,\hat{\lambda}^*)-\mathcal{L}_{\sigma}(x_{\epsilon(t)},\hat{\lambda}^*)+\frac{\epsilon(t)}{2}\Vert \hat{x}^*\Vert^2-\frac{\epsilon(t)}{2}\Vert x_{\epsilon(t)}\Vert^2\\
			\leq&\frac{\epsilon(t)}{2}\Vert \hat{x}^*\Vert^2-\frac{\epsilon(t)}{2}\Vert x_{\epsilon(t)}\Vert^2.
		\end{aligned}
	\end{equation*}
	As a consequence,
	\begin{equation*}\label{equ28}
		\frac{\epsilon(t)}{2}(\Vert x(t)-x_{\epsilon(t)}\Vert^2+\Vert x_{\epsilon(t)}\Vert^2-\Vert \hat{x}^*\Vert^2)\leq\mathcal{L}_{\epsilon(t)}(x(t))-\mathcal{L}_{\epsilon(t)}(\hat{x}^*).
	\end{equation*}
	which indicates by \eqref{equ40} that
	\begin{equation*}
		\Vert x(t)-x_{\epsilon(t)}\Vert^2+\Vert x_{\epsilon(t)}\Vert^2-\Vert \hat{x}^*\Vert^2\leq\frac{2T^2\widehat{\mathcal{G}}(T)}{t^2\beta(t)\epsilon(t)},\;\forall t\ge T.
	\end{equation*}
	Now, considering $\lim\limits_{\epsilon\rightarrow 0}\Vert x_{\epsilon}-\hat{x}^*\Vert=0$ and $\lim\limits_{t\rightarrow+\infty}t^2\beta(t)\epsilon(t)=+\infty$, it follows that
	\begin{equation*}
		\lim\limits_{t\rightarrow+\infty}\Vert x(t)-\hat{x}^*\Vert =0.
	\end{equation*}
	
	\textbf{Case II:} There exists a large enough $T$ such that $x(t)\in B(0, \Vert \hat{x}^*\Vert),\;\forall t\ge T$, i.e.,
	\begin{equation}\label{equ29}
		\Vert x(t)\Vert<\Vert \hat{x}^*\Vert,\; \forall t\ge T.
	\end{equation}
	This implies the existence of a weak cluster point for the trajectory $x(t)$, which is assumed as $\tilde{x}$ and thus there exists a sequence$(t_n)_{n\in\mathbb{N}}$ with $t_n\rightarrow+\infty$ such that $x(t_n)\rightharpoonup \tilde{x}\quad \text{as}\quad n\rightarrow+\infty$.
	Given that $f$ is continuous, the norm is weakly lower semi-continuous and $A$ is a continuous linear operator, it follows that $\mathcal{L}_{\sigma}(\cdot,\hat{\lambda}^*)$ is weakly lower semi-continuous, i.e.,
	\begin{equation*}
		\mathcal{L}_{\sigma}(\tilde{x},\hat{\lambda}^*)\leq\liminf_{n\rightarrow+\infty}\mathcal{L}_{\sigma}(x(t_n),\hat{\lambda}^*).
	\end{equation*}
	Meanwhile, from Theorem \ref{thm2.2}, we have 
	$$\lim\limits_{t\rightarrow+\infty}\mathcal{L}_{\sigma}(x(t_n),\hat{\lambda}^*)-\mathcal{L}_{\sigma}(\hat{x}^*,\hat{\lambda}^*)=0 .$$
	This indicates that 
	\begin{equation*}
		\mathcal{L}_{\sigma}(\tilde{x},\hat{\lambda}^*)\leq\mathcal{L}_{\sigma}(\hat{x}^*,\hat{\lambda}^*).
	\end{equation*}
	Thus, it follows from the fact $(\hat{x}^*,\hat{\lambda}^*)\in\Omega$ that 
	\begin{equation*}
		\mathcal{L}_{\sigma}(\tilde{x},\hat{\lambda}^*)=\mathcal{L}_{\sigma}(\hat{x}^*,\hat{\lambda}^*)=\min_{x\in\mathcal{X}}\mathcal{L}_{\sigma}(x,\hat{\lambda}^*),\quad {\rm i.e.,}\,\, 	\tilde{x}\in{\rm{argmin}}_{x\in\mathcal{X}}\mathcal{L}_{\sigma}(x,\hat{\lambda}^*)\subset S.
	\end{equation*}
	Since the norm is weakly lower semi-continuous, it follows from  \eqref{equ29} that 
	\begin{equation*}
		\Vert\tilde{x}\Vert\leq\liminf_{t\rightarrow+\infty}\Vert x(t)\Vert\leq\limsup_{t\rightarrow+\infty}\Vert x(t)\Vert\leq\Vert \hat{x}^*\Vert.
	\end{equation*}
	Keep in mind that $\hat{x}^*$ is the unique element of minimal norm in $S$ and  $\tilde{x}\in S$, and thus we can conclude that 
	\begin{equation*}
		\lim\limits_{t\rightarrow+\infty}\Vert x(t)\Vert= \Vert \tilde{x}\Vert=\Vert \hat{x}^*\Vert \quad {\rm and}  \quad	x(t)\rightharpoonup \tilde{x}= \hat{x}^*,
	\end{equation*}
	which indicates the strong convergence of the trajectory $x(t)$ i.e.,
	\begin{equation*}
		\lim\limits_{t\rightarrow+\infty}\Vert x(t)-\hat{x}^*\Vert =0.
	\end{equation*}
	
	\textbf{Case III:} For any $T\ge t_0$, there exist $T\leq s_1\neq s_2 \ge T$ such that  $x(s_1)\in B(0,\Vert \hat{x}^*\Vert)$ or $x(s_2)\notin B(0,\Vert \hat{x}^*\Vert)$. As a consequence, there exists a sequence $(t_n)_{n\in \mathbb{N}}$ such that $t_n\rightarrow+\infty$ as $n\rightarrow+\infty$, and
	\begin{equation*}
		\Vert x(t_n)\Vert =\Vert \hat{x}^*\Vert,\;\forall n\in \mathbb{N}, 
	\end{equation*}
	which can guarantee the existence of its weak sequential cluster point and without loss of generality, we assume that the whole sequence  $x(t_n)$ is weakly convergent.  With similar argument as the proof in Case II, we can obtain that
	\begin{equation*}
		\lim\limits_{n\rightarrow+\infty}\Vert x(t_n)\Vert =\Vert \hat{x}^*\Vert,\quad {\rm and}\quad  x(t_n)\rightharpoonup \hat{x}^*,\,{\rm{as}}\, n\rightarrow+\infty.
	\end{equation*}
	Therefore, 
	\begin{equation*}
		\lim\limits_{n\rightarrow+\infty}\Vert x(t_n)-\hat{x}^*\Vert =0,
	\end{equation*}
	which indicates that 
	\begin{equation*}
		\liminf_{t\rightarrow+\infty}\Vert x(t)-\hat{x}^*\Vert =0.
	\end{equation*}
	This completes the proof of Theorem \ref{thm3.1}.
	\qed
		\begin{remark}
				The $\epsilon(t)$ that satisfies {\rm{Theorem \ref{thm3.1}}} exists, for example, $\epsilon(t)=\frac{c}{t^r}\;(c>0,\;\beta<r<\beta+2)$ when $\beta(t) = t^\beta\;(\beta>0)$.
			\end{remark}

	\section{Numerical Experiments}\label{sec4}
	
	In this section, in order to illustrate the validity of the theoretical results of the proposed dynamical system \eqref{equ5}, we give the corresponding numerical experiments. All numerical experiments are run on a MacBook Pro (with Apple M2 and 16GB memory) under MATLAB version R2021b.
	
	\subsection{\bf A Convex Optimization Problem}\label{subsec7}
	In this subsection, we consider the following convex optimization problem
	\begin{equation}\label{equ42}
		\begin{aligned}
			&\min_x \; f(x) = (mx_1+nx_2+ex_3)^2,\\
			&\;{\rm{s.t.}} \; Ax =b,
		\end{aligned}
	\end{equation}
	where $f: \mathbb{R}^3\rightarrow \mathbb{R}, A = (m,-n,e), b=0$ and $m, n, e\in\mathbb{R}\backslash \left\{0\right\}$. By simple calculation, it is very easy to check that the solution set of this convex optimization problem is $\left\{(x_1, 0, -\frac{m}{e}x_1)^T\vert\; x_1\in\mathbb{R}\right\}$, $\hat{x}^*=(0,0,0)$ is the minimal norm solution of this convex optimization problem and the optimal objective function value is 0. Next, we will use the proposed dynamical system \eqref{equ5} to solve the problem.
	
\begin{figure}[H]
		\centering
		{
			\begin{minipage}[t]{0.49\linewidth}
				\centering
				\includegraphics[width=2.65in]{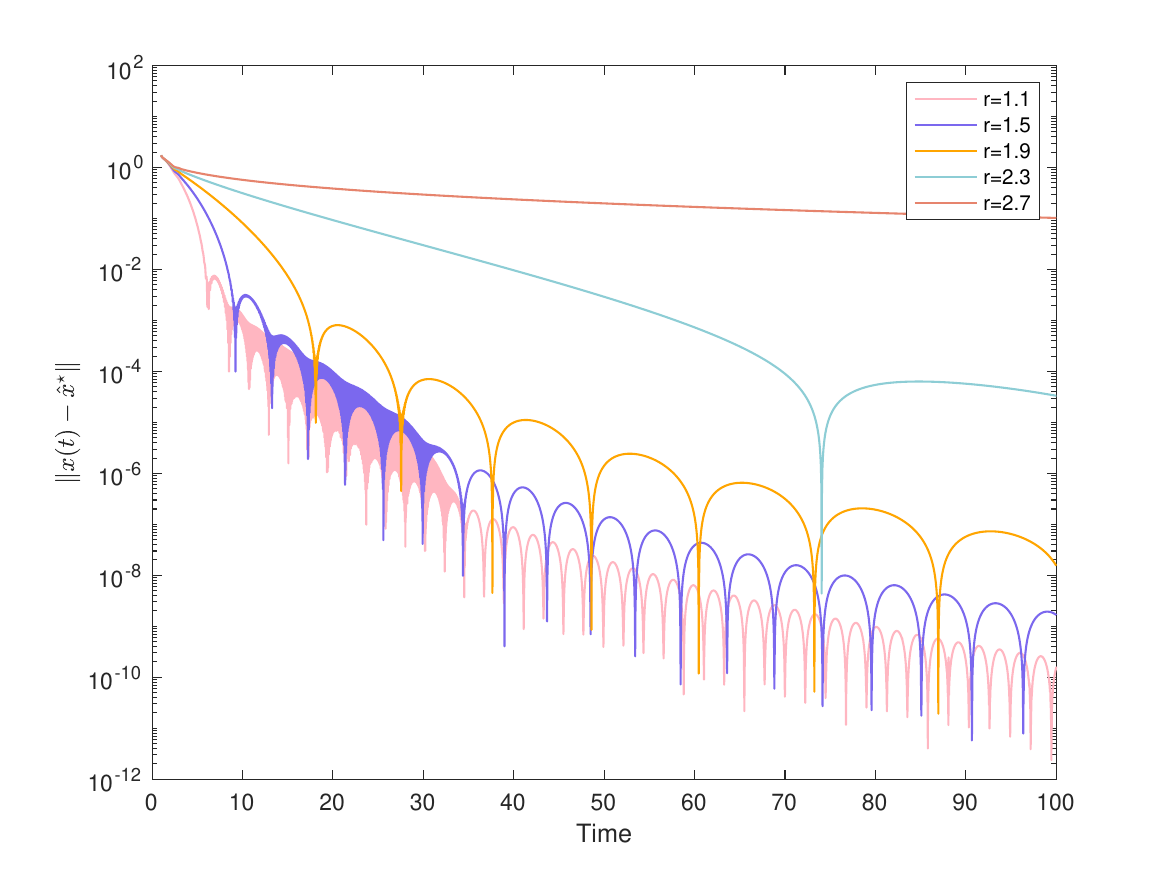}
			\end{minipage}%
		}
		{
			\begin{minipage}[t]{0.49\linewidth}
				\centering
				\includegraphics[width=2.65in]{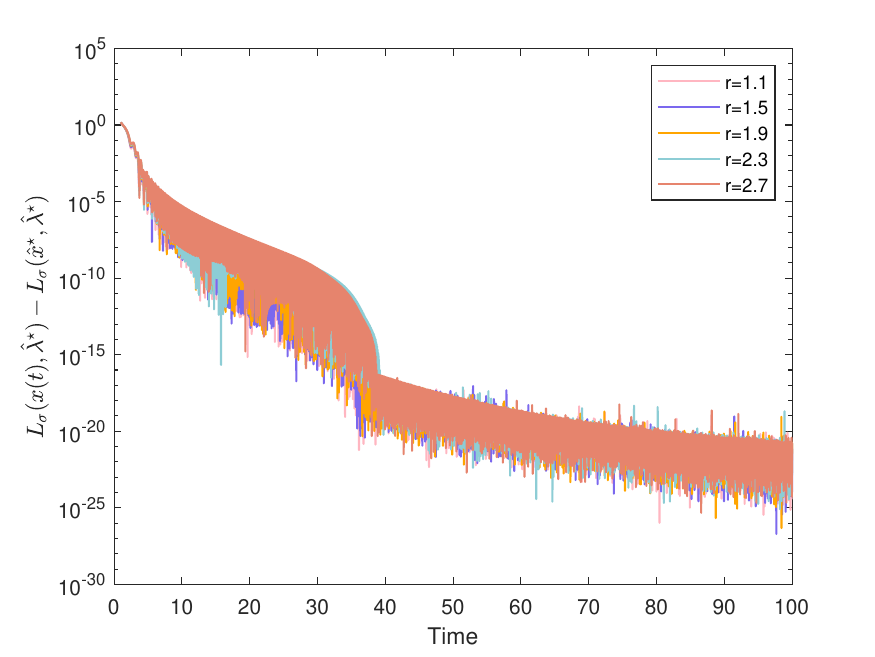}
			\end{minipage}%
		}
		\caption{Error analysis with different values of parameter $r$ in dynamical system \eqref{equ5} for problem \eqref{equ42}}
		\label{fig:1_self}
		\centering
	\end{figure}

\begin{figure}[H]
		\centering
		\subfloat[dynamical system \eqref{equ5} with Tikhonov regularization ]
		{
			\begin{minipage}[t]{0.49\linewidth}
				\centering
				\includegraphics[width=2.65in]{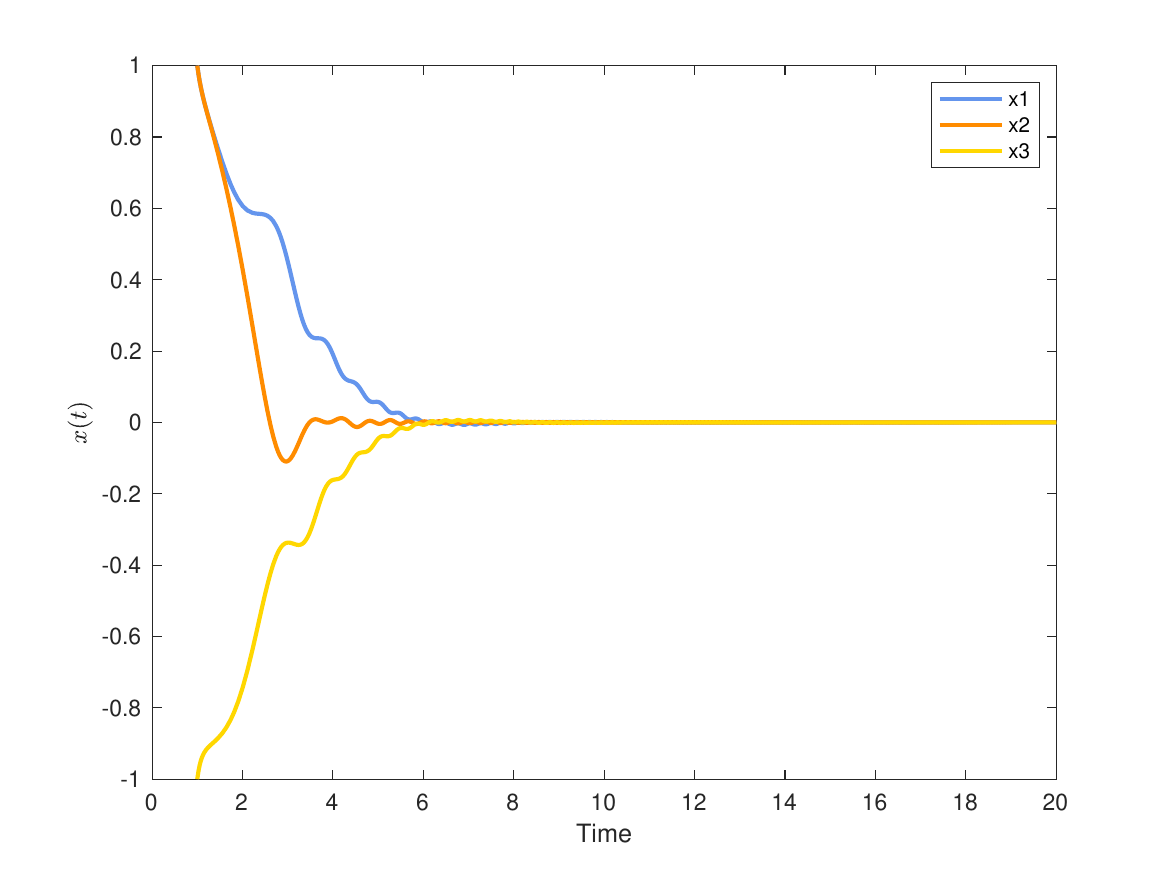}
			\end{minipage}%
		}
		\subfloat[dynamical system \eqref{equ5} without Tikhonov regularization .]
		{
			\begin{minipage}[t]{0.49\linewidth}
				\centering
				\includegraphics[width=2.65in]{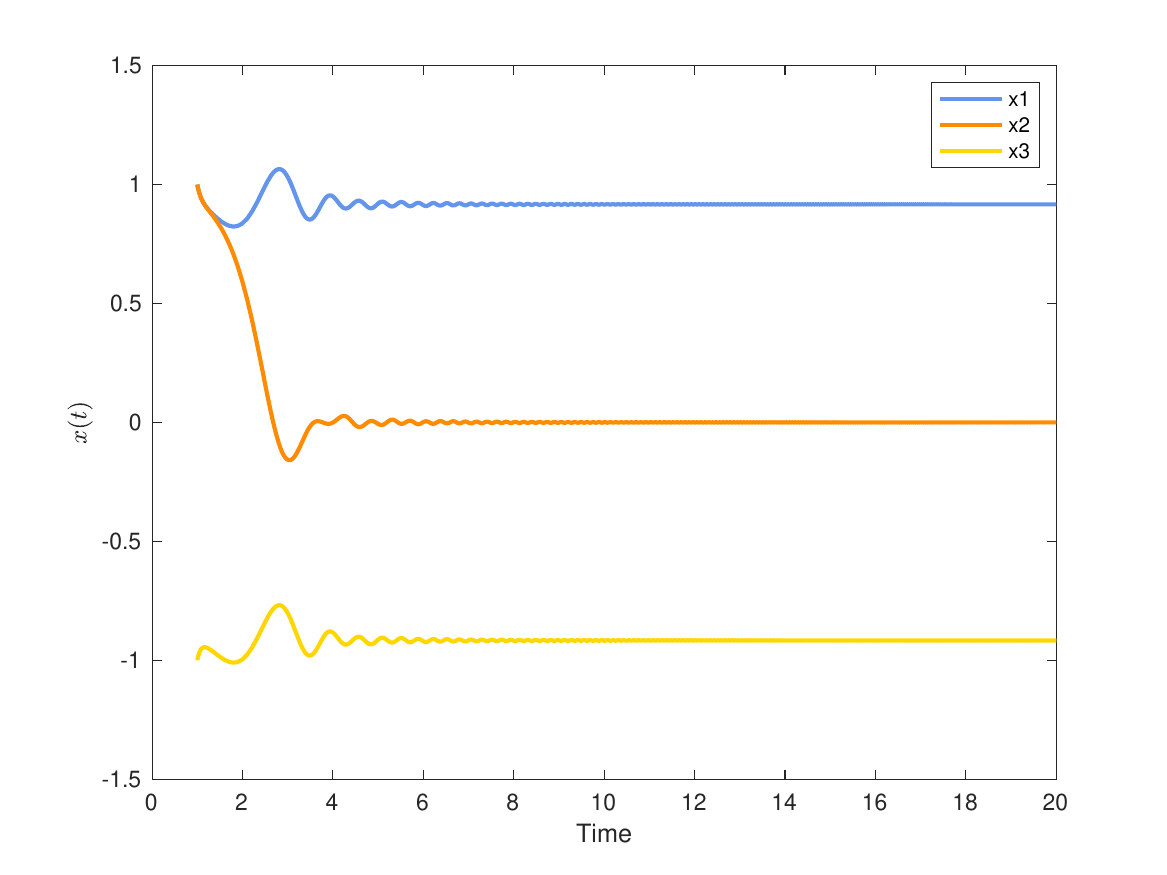}
			\end{minipage}%
		}
		\caption{The behaviors of the trajectory $x(t)$ generated by our proposed system \eqref{equ5} for problem \eqref{equ42}}
		\label{fig:1_convergence}
		\centering
	\end{figure}

	In the numerical experiments, we take $m = 1$, $n = 1$, $e = 1$, $\alpha = 13$, $\beta =1$, $\sigma=1$, $\theta=\frac{1}{\alpha-1}$, $\gamma(t) = \frac{\alpha}{t}$, $\beta(t) = t^\beta$, $\epsilon(t) = \frac{3}{t^r}$ and the start point $x(1)= (1, 1, -1)^T$, $\lambda(1) = 1$, $\dot{x}(1) = (-1, -1, 1)^T$. Next, we use the function $ode 23$ in MATLAB to solve the dynamical system \eqref{equ5}. From \cref{fig:1_self}, we can observe that $\Vert x(t)-\hat{x}^*\Vert$  along the trajectory $x(t)$ generated by the dynamical  system \eqref{equ5} varies with the parameter $r$. Specifically, as the parameter $r\;(1<r<3)$ becomes smaller, the dynamical system \eqref{equ5} has a faster descent rate and achieves higher precision in error  $\Vert x(t)-\hat{x}^*\Vert$. However, as the value of parameter $r$ changes, the variation of $ L_{\sigma}(x(t),\hat{\lambda}^*)-L_{\sigma}(\hat{x}^*,\hat{\lambda}^*) $ along the trajectory $x(t)$ generated by the dynamical system \eqref{equ5} is not obvious, that is, the error $ L_{\sigma}(x(t),\hat{\lambda}^*)-L_{\sigma}(\hat{x}^*,\hat{\lambda}^*) $ is not sensitive to the value of parameter $r$. 
	\vspace{-15pt}
	\vspace{-13pt}
	Furthermore, to analyze the behavior of the trajectory generated by the dynamical system \eqref{equ5}, we fix the parameter $r = 1.1$ in the previous experiment. 
	As illustrated in \cref{fig:1_convergence} (a), the trajectory $ x(t)$ generated by  the dynamical system \eqref{equ5} converges effectively to the minimal norm solution 
	$\hat{x}^*$. 
	In contrast, \cref{fig:1_convergence} (b) reveals that the trajectory $ x(t)$  fails to converge to the minimal norm solution $\hat{x}^*$ when the Tikhonov regularization term 
	$\epsilon(t)x(t)$ is omitted from the dynamical system \eqref{equ5}.
	\vspace{-15pt}

	Next, we conduct a comprehensive comparison between the dynamical system \eqref{equ5} and several established systems, including the Tikhonov regularized second-order plus first-order primal-dual dynamical system (TRSPFPD for short) \cite{bib23}, the dynamical system (PD-AVD) \cite{bib16}, and the dynamical system (TRIALS) \cite{bib15}.

\begin{figure}[H]
		\centering
		{
			\begin{minipage}[t]{0.49\linewidth}
				\centering
				\includegraphics[width=2.65in]{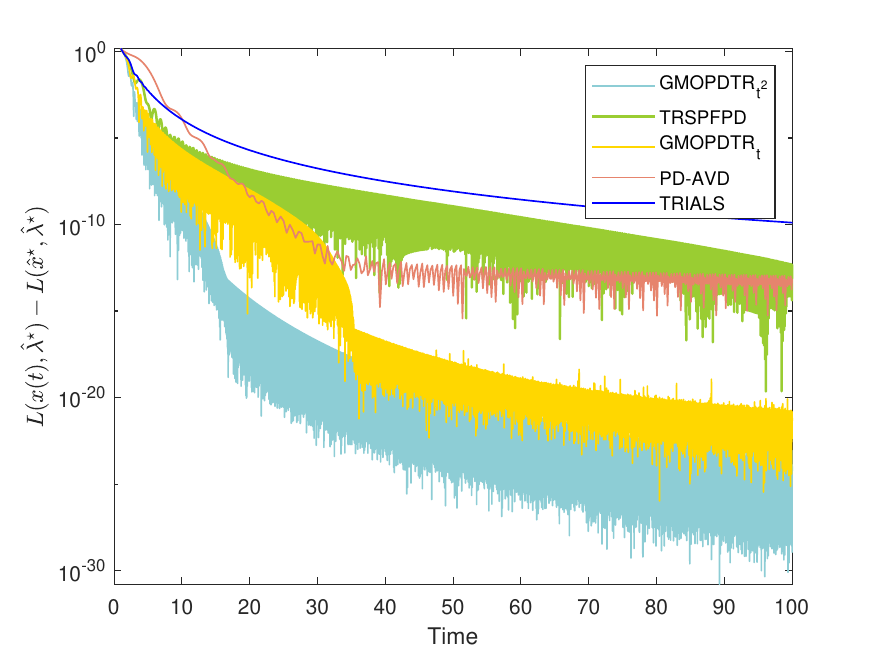}
			\end{minipage}%
		}
		{
			\begin{minipage}[t]{0.49\linewidth}
				\centering
				\includegraphics[width=2.65in]{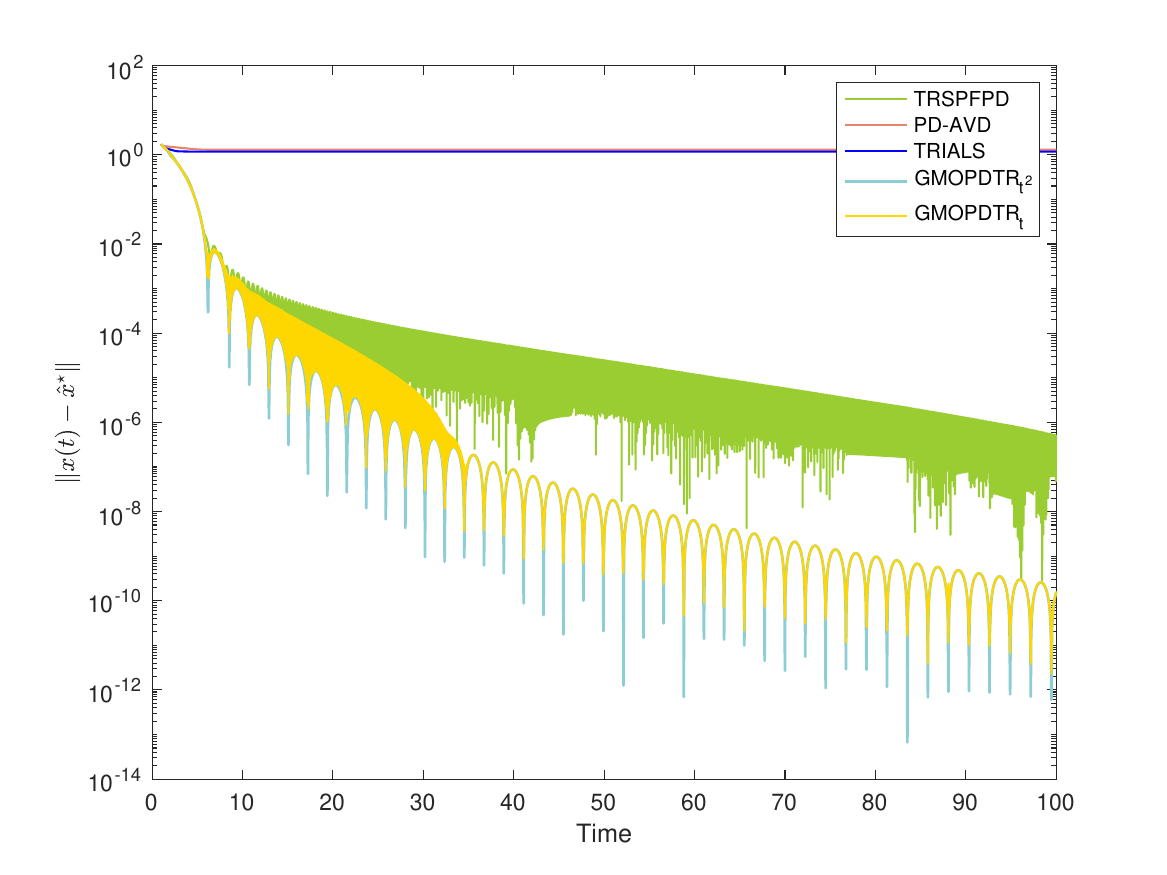}
			\end{minipage}%
		}
		\caption{Comparison of error results between the dynamical system \eqref{equ5},  (TRSPFPD), (PD-AVD) and  (TRIALS)}
		\label{fig: 1_contrast}
		\centering
	\end{figure}

	In our analysis, we set the parameters as follows: in the first case, we set 
	$m=n=e=\beta=\sigma=1$, $\alpha = 13$, $\theta=\frac{1}{\alpha-1}$, $r = 1.1$, $\gamma(t) = \frac{\alpha}{t}$, $\beta(t) = t^\beta$ and $\epsilon(t) = \frac{3}{t^r}$, while in the second case,  we maintain the same values for $m, n, e, \sigma, \alpha, \theta, \gamma(t), \beta(t)$ and $\epsilon(t)$, and we adjust $\beta$ = 2 and $r$ = 2.1. We denote the dynamical system \eqref{equ5} with these two cases as $\text{GMOPDTR}_t$ and $\text{GMOPDTR}_{t^2}$, respectively. For the dynamical system (TRSPFPD), we set 
	$\alpha =13$, $\rho = 1$,  $\epsilon(t) = \frac{3}{t^r}$ and $r = 0.1$. The initial conditions for these dynamical systems ($\text{GMOPDTR}_t$,  $\text{GMOPDTR}_{t^2}$ and TRSPFPD) are 
	$x(1)= (1, 1, -1)^T$, $\lambda(1) = 1$ and $\dot{x}(1) = (-1, -1, 1)^T$.
	For the dynamical system (PD-AVD), we set 
	$\alpha =13$, $\beta = 1$, $\theta=\frac{1}{\alpha-1}$, while for the dynamical system (TRIALS), we set 
	$\alpha_0=\frac{1}{3}$, $\eta =2$, $\alpha(t)=\alpha_0 t$, $\gamma(t)=\frac{\eta+\alpha_0}{\alpha_0t}$ and $b(t)=t^{\frac{1}{\alpha_0}-2}$. Both the dynamical system (PD-AVD) and the dynamical system (TRIALS) start with the same initial conditions: 
	$x(1)= (1, 1, -1)^T$, $\lambda(1) = 1$, $\dot{x}(1) = (-1, -1, 1)^T$ and $\dot{\lambda}(1)=-1$.  
	As shown in \cref{fig: 1_contrast}, we observe that both the error 
	$ L_{\sigma}(x(t),\hat{\lambda}^*)-L_{\sigma}(\hat{x}^*,\hat{\lambda}^*) $ and $\Vert x(t)-\hat{x}^*\Vert$ decreases more rapidly along the trajectories 
	$x(t)$ generated by the dynamical systems ($\text{GMOPDTR}_t$ and $\text{GMOPDTR}_{t^2}$) compared to those generated by the other three systems (TRSPFPD, PD-AVD and TRIALS). However, it is noteworthy that the acceleration effect in the error 
	$\Vert x(t)-\hat{x}^*\Vert$ becomes less significant in the dynamical system ($\text{GMOPDTR}_{t^{\beta}}$) when 
	$\beta\ge 1$.
	
	\subsection{\bf A Quadratic Programming Problem}\label{subsec8}
	In this subsection, we consider the quadratic programming problem defined as follows:
	\begin{equation}\label{equ43}
		\begin{aligned}
			&\min_x \; f(x) = \frac{1}{2}x^TMx+q^Tx,\\
			&\;{\rm{s.t.}} \; Ax =b,
		\end{aligned}
	\end{equation}
	where $M\in\mathbb{R}^{n\times n}$ is a symmetric and positive semi-definite matrix, $A\in\mathbb{R}^{m\times n}, q\in\mathbb{R}^n$ and $b\in\mathbb{R}^m$. 
	
\begin{figure}[H]
		\centering
		{
			\begin{minipage}[t]{0.48\linewidth}
				\centering
				\includegraphics[width=2.65in]{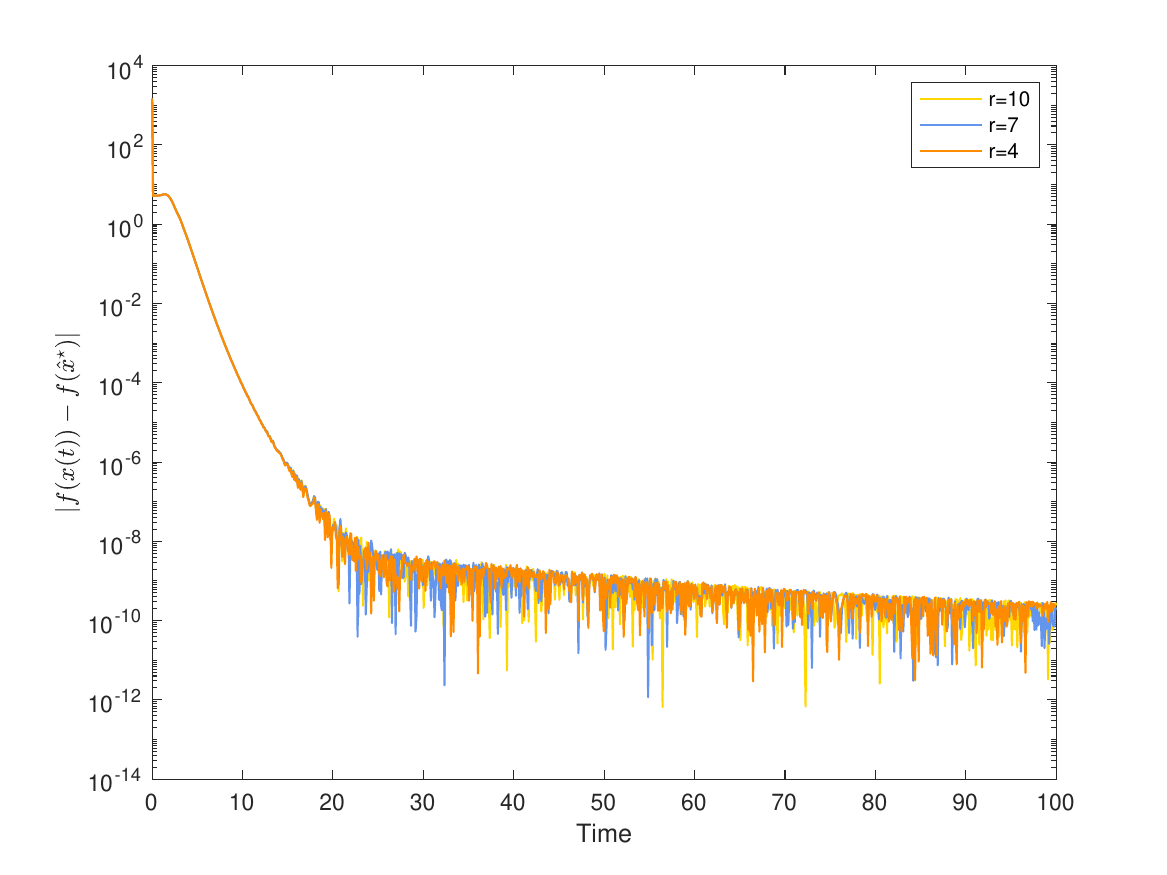}
			\end{minipage}%
		}
		{
			\begin{minipage}[t]{0.48\linewidth}
				\centering
				\includegraphics[width=2.65in]{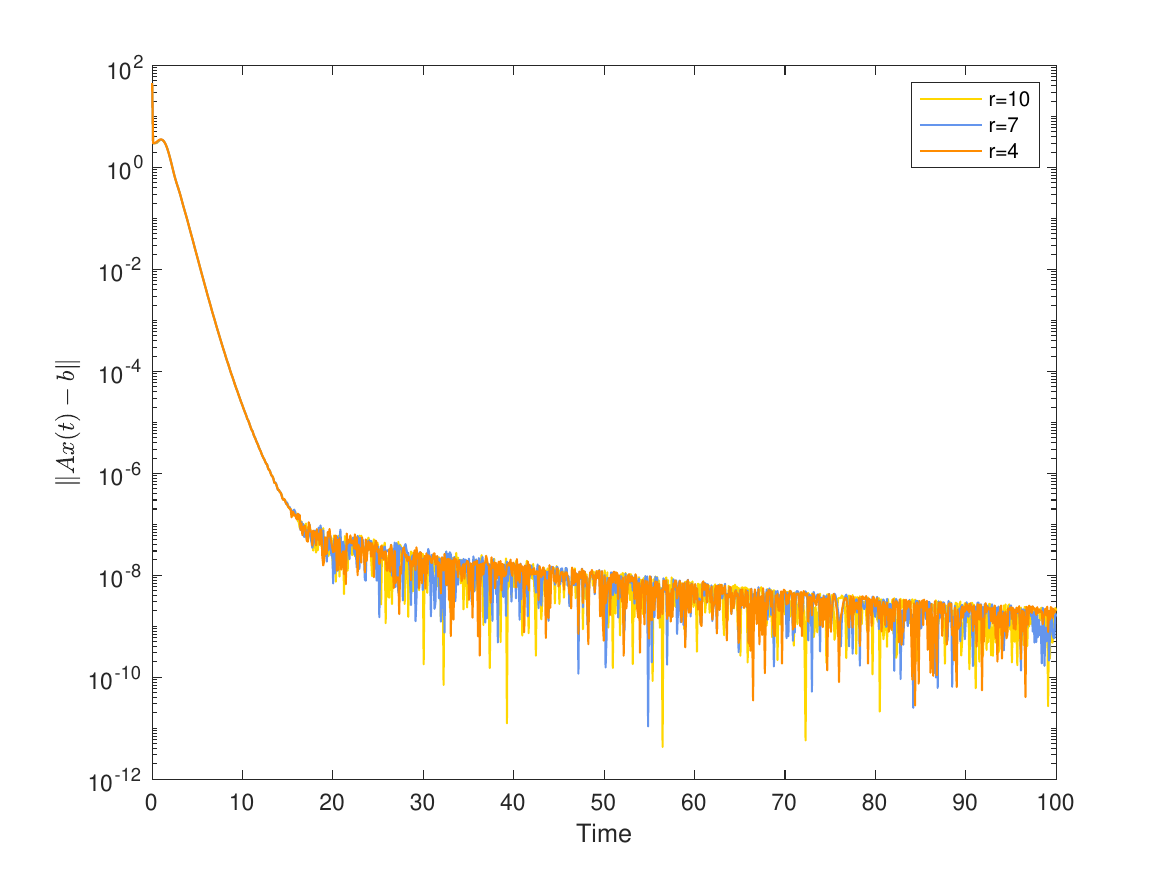}
			\end{minipage}%
		}
		\caption{Error analysis with different values of parameter $r$ in dynamical system \eqref{equ5} for problem \eqref{equ43}}
		\label{fig: 2_self}
		\centering
	\end{figure}

\begin{figure}[H]
		\centering
		{
			\begin{minipage}[t]{0.48\linewidth}
				\centering
				\includegraphics[width=2.65in]{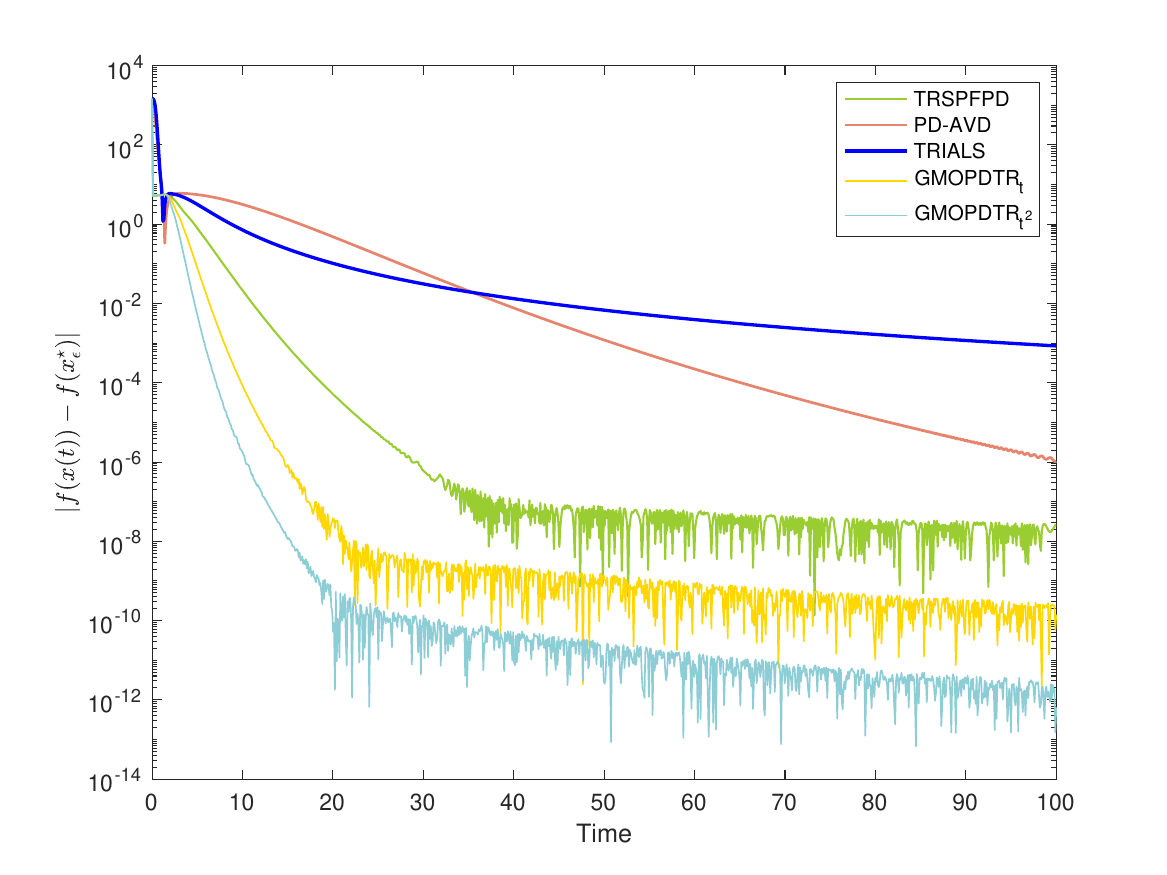}
			\end{minipage}%
		}
		{
			\begin{minipage}[t]{0.48\linewidth}
				\centering
				\includegraphics[width=2.65in]{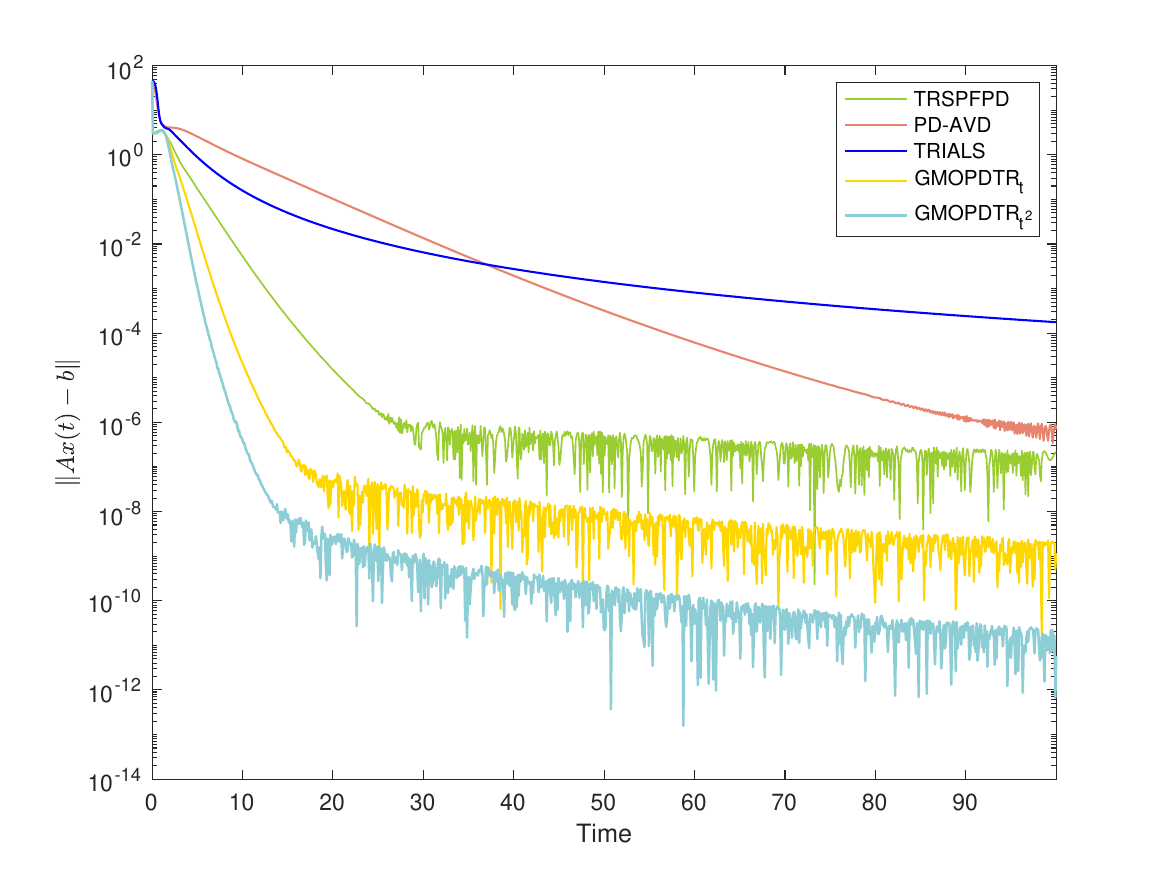}
			\end{minipage}%
		}
		\caption{Comparison of error results between the dynamical system \eqref{equ5}, (TRSPFPD), (PD-AVD) and (TRIALS)}
		\label{fig:2_contrast}
		\centering
	\end{figure}

	In the numerical experiments, we set 
$m = 30$ and $n = 50$. The vectors 
	$q$ and $A$ are generated using a standard Gaussian distribution, while the vector $b$ is generated from a uniform distribution. The matrix 
	$M$ is a symmetric and positive semi-definite matrix also created using a standard Gaussian distribution. The parameters are defined as follows: $\alpha = 13$, $\beta =\sigma=1$, $\theta=\frac{1}{\alpha-1}$, $\gamma(t) = \frac{\alpha}{t}$, $\beta(t) = t^\beta$, $\epsilon(t) = \frac{1}{t^r}$ and $r>3$. The initial conditions for the dynamical system \eqref{equ5} are given by $(x(1), \lambda(1), \dot{x}(1))=\textbf{1}^{(2n+m)\times 1}$. 
	Subsequently, we utilize the function $ode 23$ in MATLAB to solve the dynamical system \eqref{equ5}. Additionally, we employ the function $quadprog$ in MATLAB to derive the optimal value $f(\hat{x}^*)$ for the problem \eqref{equ43}. 
	As illustrated in \cref{fig: 2_self}, the changes in the parameter 
	$r$ do not cause significant variations in the values of
$\vert f(x(t))-f(\hat{x}^*) \vert$ and $\Vert Ax(t)-b\Vert$ along the trajectory generated by the dynamical system \eqref{equ5} . This observation indicates that the errors 
$\vert f(x(t))-f(\hat{x}^*) \vert$ and $\Vert Ax(t)-b\Vert$ exhibit limited sensitivity to variations in the parameter $r$.
	Besides, we compare the dynamical system \eqref{equ5} with TRSPFPD, PD-AVD and TRIALS. We set the parameters for the dynamical system \eqref{equ5}  as follows: for the dynamical system ($\text{GMOPDTR}_t$), we set $\alpha = 13$, $\beta =1$, $\sigma=1$, $\theta=\frac{1}{\alpha-1}$, $\gamma(t) = \frac{\alpha}{t}$, $\beta(t) = t^\beta$, $\epsilon(t) = \frac{1}{t^r}$ and $r = 4$, while for the dynamical system ($\text{GMOPDTR}_{t^2}$), we set $\alpha = 13$, $\beta =2$, $\sigma=1$, $\theta=\frac{1}{\alpha-1}$, $\gamma(t) = \frac{\alpha}{t}$, $\beta(t) = t^\beta$, $\epsilon(t) = \frac{1}{t^r}$ and $r = 5$.  
	For the dynamical system (TRSPFPD), we set $\alpha =13$, $\rho = 1$,  $\epsilon(t) = \frac{1}{t^r}$ and $r = 4$. For the dynamical system (PD-AVD), we set $\alpha =13$, $\beta = 1$ and $\theta=\frac{1}{\alpha-1}$, while for the dynamical system (TRIALS), we set $\alpha_0=\frac{1}{3}$, $\eta =2$, $\alpha(t)=\alpha_0 t$, $\gamma(t)=\frac{\eta+\alpha_0}{\alpha_0t}$and $b(t)=t^{\frac{1}{\alpha_0}-2}$. The initial conditions for these dynamical systems ($\text{GMOPDTR}_t$, $\text{GMOPDTR}_{t^2}$ and TRSPFPD) are $(x(1), \lambda(1), \dot{x}(1))=\textbf{1}^{(2n+m)\times 1}$, while the initial conditions for  the dynamical system (PD-AVD) and  the dynamical system (TRIALS) are $(x(1), \lambda(1), \dot{x}(1), \dot{\lambda}(1))=\textbf{1}^{2(n+m)\times 1}$.
	We then employ the function $ode 23$ in MATLAB to solve each dynamical system. As shown in \cref{fig:2_contrast}, both errors
$ L_{\sigma}(x(t),\hat{\lambda}^*)-L_{\sigma}(\hat{x}^*,\hat{\lambda}^*) $ and $\Vert Ax(t)-b\Vert$  decrease more rapidly along the trajectories generated by the dynamical systems ($\text{GMOPDTR}_t$ and $\text{GMOPDTR}_{t^2}$) compared to those generated by the other three dynamical systems. Additionally, in terms of computational accuracy, the dynamical systems ($\text{GMOPDTR}_t$ and $\text{GMOPDTR}_{t^2}$) outperforms the other three dynamical systems.

	\section{Conclusion}\label{sec5}

	This paper investigates a general mixed-order primal-dual dynamical system with Tikhonov regularization for a convex optimization problem with linear equality constraints. We obtain the following asymptotic convergence properties by constructing appropriate Lyapunov functions:
	\begin{itemize}
		\item[(1)] When the Tikhonov regularization parameter $\epsilon(t)$  satisfies the condition 
		$\int_{t_0}^{+\infty}t\beta(t)\epsilon(t)<+\infty$,  we obtain a convergence rate of $\mathcal{O}(\frac{1}{t^2\beta(t)})$ for the primal-dual gap, feasibility violation, and objective function error along the trajectory generated by the dynamical system \eqref{equ5}.
		\item[(2)] When the Tikhonov regularization parameter $\epsilon(t)$  satisfies the condition $\int_{t_0}^{+\infty}\frac{\beta(t)\epsilon(t)}{t}<+\infty$, we demonstrate that the primal-dual gap converges at a rate of 
		$o(\frac{1}{\beta(t)})$. Furthermore, the trajectory generated by the dynamical system \eqref{equ5} strongly converges to the minimum norm solution of the optimization problem \eqref{equ1}.
	\end{itemize}
	We conclude our study with numerical experiments that validate our theoretical findings. These experiments compare the performance of our dynamical system \eqref{equ5} with those of the dynamical systems TRSPFPD proposed by Zhu et al. \cite{bib23}, PD-AVD introduced by Bot et al. \cite{bib16}, and TRIALS proposed by Attouch et al. \cite{bib15}, all addressing the same optimization problem. The results of these experiments reveal that our dynamical system \eqref{equ5} not only exhibits a faster convergence rate but also attains higher solution accuracy, particularly evident in the second experiment.

	\begin{acknowledgements}
		This research was supported by the National Natural Science Foundation of China (12171070, 12401401) and Sichuan Science and Technology Program (Grant No. 2025ZNSFSC0813).
	\end{acknowledgements}
	
	\begin{appendices}
		
		\section{Well-Posedness of the Cauchy Problem}\label{Appendix A}
		
		In this appendix, we will give the proof of the existence and uniqueness of a global solution for dynamical system \eqref{equ5} with the initial condition $x(t_0)=x_0$, $\lambda(t_0)=\lambda_0$, $\dot{x}(t_0)=u_0$.
		\vskip 2mm
		\begin{definition}
			A function $x : [t_0, +\infty)\rightarrow\mathcal{X}$ is a strong global solution of dynamical system \eqref{equ5} if it satisfies the following properties:
			\begin{enumerate}
				\item[(i)] $x : [t_0, +\infty)\rightarrow\mathcal{X}$ and $\lambda : [t_0, +\infty)\rightarrow\mathcal{Y}$ are locally absolutely continuous;
				\item[(ii)] 		
				$
				\begin{cases}
					\ddot{x}(t) + \gamma(t) \dot{x}(t) + \beta(t) (\nabla f(x(t)) + A^T \lambda(t) + \sigma A^T (Ax(t) - b) + \epsilon(t) x(t)) = 0, \\
					\dot{\lambda}(t) - t \beta(t) (A(x(t) + \theta t \dot{x}(t)) - b) = 0;
				\end{cases}
				,\forall t\in[t_0, +\infty);
				$
				\item[(iii)] $x(t_0) = x_0\in\mathcal{X}$, $\lambda(t_0) = \lambda_0\in\mathcal{Y}$,  $\dot{x}(t_0) = u_0\in\mathcal{X}$.
			\end{enumerate}
			\vskip 2mm
		\end{definition}
		\begin{theorem}\label{thm1}
			Suppose that $f$ is a continuous function such that $\nabla f$ is Lipschitz continuous with Lipschitz constant $L>0$. Let $\epsilon,\gamma,\beta:\left[t_0,+\infty\right)\rightarrow\left(0,+\infty\right)$ be continuous functions. Then, for any initial condition  $(x(t_0),\lambda(t_0),\dot{x}(t_0))=(x_0,\lambda_0,u_0)\in \mathcal{X}\times\mathcal{Y}\times\mathcal{X}$, the dynamical system \eqref{equ5} has a unique strong global solution.
		\end{theorem}
		{\it Proof} 
		Denote $Y(t)$ := $(y_1(t),y_2(t),y_3(t))$ = $(x(t),\lambda(t),\dot{x}(t))$, $Y_0$ := $(x_0,\lambda_0,u_0)$. For the ease of writing, we denote $Y(t)$ as $Y$ := $(y_1,y_2,y_3)$. Then, the dynamical system \eqref{equ5} can be rewritten as: 
		\begin{equation}\label{equA1}
			\begin{cases}
				\frac{dY}{dt}=F(t,Y),\\
				Y(t_0)=Y_0,
			\end{cases}
		\end{equation}
		where
		\begin{equation*}
			F(t,Y) := \left(
			\begin{aligned}
				&y_3\\
				&t\beta(t)(A(y_1+\theta ty_3)-b)\\
				&-\gamma(t)y_3-\beta(t)(\nabla f(y_1)+A^Ty_2+\sigma A^T(Ay_1-b)+\epsilon(t)y_1)
			\end{aligned}\right).
		\end{equation*}
		Since $\nabla f$ is Lipschitz continuous and $A$ is linear, it follows that for any $Y$, $\widetilde{Y}\in \mathcal{X}\times\mathcal{Y}\times\mathcal{X}$,
		\begin{eqnarray*}
			\Vert F(s,Y(s))-F(s,\widetilde{Y}(s))\Vert&\leq&(1+\gamma(t)+\theta t^2\beta(t)\Vert A\Vert) \Vert y_3-\tilde{y}_3\Vert+\beta(t)\Vert A\Vert \Vert \tilde{y}_2-y_2\Vert\\
			&&+(\sigma\beta(t)\Vert A^TA\Vert+\beta(t)\epsilon(t)+t\beta(t)\Vert A\Vert) \Vert\tilde{y}_1-y_1\Vert\\
			&&+\beta(t)\Vert \nabla f(\tilde{y}_1)-\nabla f(y_1)\Vert\\
			&\leq&(1+\gamma(t)+\theta t^2\beta(t)\Vert A\Vert) \Vert y_3-\tilde{y}_3\Vert+\beta(t)\Vert A\Vert \Vert \tilde{y}_2-y_2\Vert\\
			&&+(\sigma\beta(t)\Vert A^TA\Vert+\beta(t)\epsilon(t)+t\beta(t)\Vert A\Vert+L\beta(t)) \Vert\tilde{y}_1-y_1\Vert\\
			&\leq&M(t)\Vert Y-\widetilde{Y}\Vert,
		\end{eqnarray*}
		where $C$ := $L+\sigma\Vert A^TA\Vert+\Vert A\Vert$ and $M(t)$ := $1+C\beta(t)+\gamma(t)+\Vert A\Vert \theta t^2\beta(t)+\Vert A\Vert t\beta(t)+\beta(t)\epsilon(t)$.
		Since $\epsilon,\gamma,\beta:\left[t_0,+\infty\right)\rightarrow\left(0,+\infty\right)$ are continuous functions, we have $M(t)\in L_{loc}^{1}\left[t_0,+\infty\right)$. Further, by the Lipschitz continuity of 
		$\nabla f$, we deduce that
		\begin{equation*}
			\Vert \nabla f(y_1)\Vert\leq\Vert\nabla f(0)\Vert+ L\Vert y_1\Vert.
		\end{equation*}
		Therefore, for any given $Y\in\mathcal{X}\times\mathcal{Y}\times\mathcal{X}$ and $t_0< T <+\infty$, the following inequality holds:
		\begin{eqnarray*}
			\int_{t_0}^{T}\Vert F(t,Y)\Vert dt&\leq&\int_{t_0}^{T}(\sigma\beta(t)\Vert A^TA\Vert+\beta(t)\epsilon(t)+t\beta(t)\Vert A\Vert+L\beta(t))\Vert y_1\Vert+\nabla f(0)\beta(t)+\sigma\Vert A^Tb\Vert\beta(t)\\
			&&\quad+\beta(t)\Vert A\Vert\Vert y_2\Vert +(1+\gamma(t)+\theta t^2\beta(t)\Vert A\Vert)\Vert y_3\Vert dt\\
			&\leq&\int_{t_0}^{T}(\sigma\beta(t)\Vert A^TA\Vert+\beta(t)\epsilon(t)+t\beta(t)\Vert A\Vert+L\beta(t)+\nabla f(0)\beta(t)+\sigma\Vert A^Tb\Vert\beta(t))(1+\Vert Y\Vert)\\
			&&\quad+\beta(t)\Vert A\Vert(1+\Vert Y\Vert)+(1+\gamma(t)+\theta t^\beta(t)\Vert A\Vert)(1+\Vert Y\Vert)dt\\
			&\leq& \int_{t_0}^{T}S(t)(1+\Vert Y\Vert)dt,
		\end{eqnarray*}
		where $S(t)$ := $1+C\beta(t)+\gamma(t)+\Vert A\Vert \theta t^2\beta(t)+\Vert A\Vert t\beta(t)+\beta(t)\epsilon(t)$ with $C$ := $L+\sigma\Vert A^TA\Vert+\Vert A\Vert+\nabla f(0)+\sigma\Vert A^Tb\Vert$. It is clear that $S(t)\in L_{loc}^{1}\left[t_0,+\infty\right)$ and for any $Y\in\mathcal{X}\times\mathcal{Y}\times\mathcal{X}$, $F(\cdot, Y)\in L_{loc}^{1}([t_0,+\infty); \mathcal{X}\times\mathcal{Y}\times\mathcal{X})$. According to \cite[Proposition 6.2.1]{bib24} and \cite[Theorem 5]{bib15}, the  Cauchy problem \eqref{equA1} has a unique global solution $Y\in W_{loc}^{1,1}([t_0,+\infty); \mathcal{X}\times\mathcal{Y}\times\mathcal{X})$. Therefore, the dynamical system \eqref{equ5} has a unique strong global solution $(x, \lambda)$. This completes the proof of Theorem \ref{thm1}.
		\qed

		\section{Some Auxiliary Results}\label{secA1}
		
		\begin{lemma}{\rm{\cite[Lemma 6]{bib17}}}\label{lemB.1}
			Assume that $g:\left[t_0,+\infty\right)\rightarrow\mathcal{Y}$ is a continuous differentiable function, $a:\left[t_0,+\infty\right)\rightarrow\left[0,+\infty\right)$ is also continuous differentiable, $t_0>0$ and $C\ge 0$. If
			\begin{equation*}
				\Vert g(t)+\int_{t_0}^{t} a(s)g(s) ds\Vert \leq C,\quad\forall t\ge t_0,
			\end{equation*}
			then
			\begin{equation*}
				\sup_{t\ge t_0}\Vert g(t)\Vert<+\infty.
			\end{equation*}
		\end{lemma}
		
		\begin{lemma}{\rm{\cite[Lemma A.3]{bib21}}}\label{lemB.2}
			Suppose that $\delta>0$ and $\phi\in L^1\left(\delta,+\infty\right)$ is a non-negative and continuous function. Additionally, let $\psi:\left[\delta,+\infty\right)\rightarrow\left(0,+\infty\right)$ be a non-decreasing function satisfying $\lim\limits_{t\rightarrow+\infty}\psi(t)=+\infty$. Then,
			\begin{equation*}
				\lim\limits_{t\rightarrow+\infty}\frac{1}{\psi(t)}\int_{\delta}^{t}\psi(s)\phi(s)ds =0.
			\end{equation*}
		\end{lemma}
		
		\begin{lemma}\label{lemB.3}
			Let $\theta>0$ and $z\in \mathcal{X}$ be a fixed element, and let $x:\left[t_0,+\infty\right)\rightarrow\mathcal{X}$ be a continuously differentiable function. If there exists a constant $\widehat{C}$ such that
			\begin{equation*}
				\Vert x(t)-z+\theta t\dot{x}(t)\Vert^2\leq\widehat{C},\quad\forall t\in\left[t_0,+\infty\right),
			\end{equation*}
			then $x(t)$ is bounded on $\left[t_0,+\infty\right)$.
		\end{lemma}
		{\it Proof} 
		Under the given conditions, the following inequality holds:
		\begin{equation*}
			\Vert x(t)-z\Vert^2+2\theta t\langle x(t)-z,\dot{x}(t)\rangle\leq \widehat{C}.
		\end{equation*}
		Dividing both sides of the above inequality by $p(t)$ := $(\frac{t}{t_0})^{\frac{\theta+1}{\theta}}$ yields
		\begin{equation*}
			\frac{\Vert x(t)-z\Vert^2}{p(t)}+\frac{2\theta t}{p(t)}\langle x(t)-z,\dot{x}(t)\rangle\leq\frac{\widehat{C}}{p(t)}.
		\end{equation*}
		Let $h(t)$ := $\frac{1}{2}\Vert x(t)-z\Vert^2$. It follows that $\dot{h}(t)=\langle x(t)-z,\dot{x}(t)\rangle$. Substituting them into the previous inequality leads to
		\begin{equation}\label{equB2}
			\frac{h(t)}{p(t)}+q(t)\dot{h}(t)\leq\frac{\widehat{C_1}}{p(t)},
		\end{equation}
		where $q(t)$ := $\frac{\theta t}{p(t)}$ and $\widehat{C_1}$ := $\frac{\widehat{C}}{2}$. From the definition of $p(t)$, it follows that $\dot{q}(t)=-\frac{1}{p(t)}$, indicating that $q(t)$ is bounded. Rewriting  \eqref{equB2} gives
		\begin{equation*}
			q(t)\dot{h}(t)-\dot{q}(t)(h(t)-\widehat{C_1})\leq 0.      
		\end{equation*}
		Dividing the inequality by $q(t)^2$ and integrating the obtained inequality from $t_0$ to $t$ result in
		\begin{equation*}
			h(t)\leq\frac{h(t_0)-\widehat{C_1}}{q(t_0)}q(t)+\widehat{C_1}.
		\end{equation*}
		Based on the boundedness of $q(t)$ and the definition of $h(t)$ , we conclude that  $x(t)$ is bounded. This completes the proof of Lemma \ref{lemB.3}.
		\qed

	\end{appendices}


\begin{thebibliography}{}
		\bibitem{bib3}
		Alvarez, F.: On the minimizing property of a second order dissipative system in Hilbert spaces. SIAM J. Control Optim. {\bf 38}(4), 1102--1119  (2000)
		
		\bibitem{bib11}
		Aujol, J.F., Dossal, C., Rondepierre, A.: Optimal convergence rates for Nesterov acceleration. SIAM J. Optim. {\bf 29}(4), 3131--3153 (2019)
		
		\bibitem{bib50}
		Aujol, J.F., Dossal, C., Ho{\`a}ng, V.H., Labarri{\`e}re, H., Rondepierre, A.: Fast convergence of inertial dynamics with Hessian-driven damping under geometry assumptions. Appl. Math. Optim. {\bf 88}(3), 81 (2023)
		
		\bibitem{bib37}
		Alecsa, C.D., L{\'a}szl{\'o}, S.C.: Tikhonov regularization of a perturbed heavy ball system with vanishing damping. SIAM J. Optim. {\bf 31}, 2921--2954 (2021)
		
		\bibitem{bib48}
		Alecsa, C.D., L{\'a}szl{\'o}, S.C., Pinta, T.: An extension of the second order dynamical system that models Nesterov’s convex gradient method. Appl. Math. Optim. {\bf 84}, 1687--1716 (2021)
		
		\bibitem{bib9}
		Attouch, H., Chbani, Z., Riahi, H.: Rate of convergence of the Nesterov accelerated gradient method in the subcritical case $\alpha\leq$3. ESAIM Control Optim. Calc. Var. {\bf 25}, 2 (2019)
		
		\bibitem{bib12}
		Attouch, H., Chbani, Z., Riahi, H.: Fast convex optimization via time scaling of damped inertial gradient dynamics. Pure Appl. Funct. Anal. {\bf 6}(6), 1081--1117 (2021)
		
		\bibitem{bib15}
		Attouch, H., Chbani, Z., Fadili, J., Riahi, H.: Fast convergence of dynamical ADMM via time scaling of damped inertial dynamics. J. Optim. Theory Appl. {\bf 193}(1-3), 704--736 (2022)
		
		\bibitem{bib19}
		Attouch, H.: Viscosity solutions of minimization problems. SIAM J. Optim. {\bf 6}(3), 769--806 (1996)
		
		\bibitem{bib20}
		Attouch, H., Czarnecki, M.-O.: Asymptotic control and stabilization of nonlinear oscillators with non-isolated equilibria. J. Differ. Equ. {\bf 179}(1), 278--310 (2002)
		
		\bibitem{bib21}
		Attouch, H., Chbani, Z., Riahi, H.: Combining fast inertial dynamics for convex optimization with Tikhonov regularization. J. Math. Anal. Appl. {\bf 457}(2), 1065--1094 (2018)
		
		\bibitem{bib25}
		Attouch, H., Cominetti, R.: A dynamical approach to convex minimization coupling approximation with the steepest descent method. J. Differ. Equ. {\bf 128}(2), 519--540 (1996)
		
		\bibitem{bib38}
		Attouch, H., Balhag, A., Chbani, Z., Riahi, H.: Damped inertial dynamics with vanishing Tikhonov regularization: strong asymptotic convergence towards the minimum norm solution. J. Differ. Equ. {\bf 311}, 29--58 (2022)
		
		\bibitem{bib39}
		Attouch, H., L{\'a}szl{\'o}, S.:  Convex optimization via inertial algorithms with vanishing Tikhonov regularization: fast convergence to the minimum norm solution. Math. Meth. Oper. Res. {\bf 99}(3), 307--347 (2024)
		
		\bibitem{bib41}
		Attouch, H., Cabot, A.: Asymptotic stabilization of inertial gradient dynamics with time-dependent viscosity. J. Differ. Equ. {\bf 263}(9), 5412--5458 (2017)
		
		\bibitem{bib42}
		Attouch, H., Cabot, A., Chbani, Z., Riahi, H.: Rate of convergence of inertial gradient dynamics with time-dependent viscous damping coefficient. Evol. Equ. Control Theory {\bf 7}(3), 353--371 (2018)
		
		\bibitem{bib43}
		Attouch, H., Chbani, Z., Peypouquet, J., Redont, P.: Fast convergence of inertial dynamics and algorithms with asymptotic vanishing viscosity. Math. Program. {\bf 168}, 123--175 (2018)
		
		\bibitem{bib49}
		Attouch, H., Balhag, A., Chbani, Z., Riahi, H.: Accelerated gradient methods combining Tikhonov regularization with geometric damping driven by the Hessian. Appl. Math. Optim. {\bf 88}(2), 29 (2023)
		
		
		\bibitem{bib51}
		Attouch, H., Peypouquet, J., Redont, P.: Fast convex optimization via inertial dynamics with Hessian driven damping. J. Differ. Equ. {\bf 261}(10), 5734--5783 (2016)
		\bibitem{bib4}
		B{\'e}gout, P., Bolte, J., Jendoubi, M.A.: On damped second-order gradient systems. J. Differ. Equ. {\bf 259}(7), 3115--3143 (2015)
		
		\bibitem{bib26}
		Boyd, S., Parikh, N., Chu, E., Peleato, B., Eckstein, J.: Distributed optimization and statistical learning via the alternating direction method of multipliers. Found. Trends Mach. Learn. {\bf 3}, 1--122 (2011)
		
		\bibitem{bib33}
		Browder, F.E.: Existence and approximation of solutions of nonlinear variational inequalities. Proc. Nutl. Acad. Sci. {\bf 56}(4), 108--1086 (1966)
		
		\bibitem{bib58}
		Botsaris, C.A., Jacobson, D.H.: A Newton-type curvilinear search method for optimization. Math. Anal. Appl. {\bf 54}(1), 217--229 (1976)
		
		\bibitem{bib16}
		Bo{\c{t}}, R.I., Nguyen, D.-K.: Improved convergence rates and trajectory convergence for primal-dual dynamical systems with vanishing damping. J. Differ. Equ. {\bf 303}, 369--406 (2021)
		
		
		\bibitem{bib40}
		Bo{\c{t}}, R.I., Csetnek, E.R., László, S.C.: Tikhonov regularization of a second order dynamical system with Hessian driven damping. Math. Program. {\bf 189}, 151--186 (2021)
		
		\bibitem{bib18}
		Cominetti, R., Peypouquet, J., Sorin, S.: Strong asymptotic convergence of evolution equations governed by maximal monotone operators with Tikhonov regularization. J. Differ. Equ. {\bf 245}(12), 3753--3763 (2008)
		
		\bibitem{bib35}
		Cabot, A.: Proximal point algorithm controlled by a slowly vanishing term: applications to hierarchical minimization. SIAM J. Optim. {\bf 15}, 555--572 (2005)
		
		\bibitem{bib59}
		Cabot, A., Engler, H., Gadat, S.: On the long time behavior of second order differential equations with asymptotically small dissipation. Trans. Amer. Math. Soc. {\bf 361}(11), 5983--6017 (2009)
		
		\bibitem{bib55}
		Cai, Z., Wang, Y.: A multiobjective optimization-based evolutionary algorithm for constrained optimization. IEEE Trans.  Evol. Comput. {\bf 10}(6), 658--675 (2006)
		
		
		\bibitem{bib57}
		Eriskin, L., Karatas, M., Zheng, Y.-J.: A robust multi-objective model for healthcare resource management and location planning during pandemics. Ann. Oper. Res. {\bf 335}(3), 1471--1518 (2024)
		
		\bibitem{bib28}
		Goldstein, T., O'Donoghue, B., Setzer, S., Baraniuk, R.: Fast alternating direction optimization methods. SIAM J. Imaging Sci. {\bf 7}, 1588--1623 (2014)
		
		\bibitem{bib24}
		Haraux, A.: Systemes dynamiques dissipatifs et applications. Recherches Math. Appl. {\bf 17}, Masson, Paris (1991)
		
		\bibitem{bib14}
		He, X., Hu, R., Fang, Y.P.: Convergence rates of inertial primal-dual dynamical methods for separable convex optimization problems. SIAM J. Control Optim. {\bf 59}(5), 3278--3301 (2021)
		
		
		\bibitem{HeTAC2022}
		He, X., Hu, R., Fang, Y.P.: ``Second-order primal'' + ``first-order dual'' dynamical systems with time scaling for linear equality constrained convex optimization problems. IEEE Trans. Automat. Control {\bf 67}(8), 4377--4383 (2022)
		
		\bibitem{bib17}
		He, X., Hu, R., Fang, Y.P.: Fast primal-dual algorithm via dynamical system for a linearly constrained convex optimization problem. Automatica {\bf 146}, 110547 (2022)
		
		
		\bibitem{bib47}
		Luo, H., Chen, L.: From differential equation solvers to accelerated first-order methods for convex optimization. Math. Program. {\bf 195}, 735--781 (2022)
		
		
		\bibitem{bib36}
		L{\'a}szl{\'o}, S.C.: On the strong convergence of the trajectories of a Tikhonov regularized second order dynamical system with asymptotically vanishing damping. J. Differ. Equ. {\bf 362}, 355--381 (2023) 
		
		\bibitem{bib8}
		May, R.: Asymptotic for a second-order evolution equation with convex potential and vanishing damping term. Turkish J. Math. {\bf 41}(3), 681--685 (2017)
		
		
		\bibitem{bib1}
		Polyak, B.T.: Some methods of speeding up the convergence of iteration methods. USSR Comput. Math. Math. Phys. {\bf 4}, 1--17 (1964)
		
		\bibitem{bib7}
		Su, W., Boyd, S., Candes, E.: A differential equation for modeling Nesterov’s accelerated gradient method: theory and insights. J. Mach. Learn. Res. {\bf 17}(153), 1--43 (2016)
		
		
		\bibitem{bib29}
		Shi, G., Johansson, K.H.: Randomized optimal consensus of multi-agent systems. Automatica {\bf 48}, 3018--3030 (2012)
		
		\bibitem{bib46}
		Shi, B., Du, S.S., Jordan, M.I., Su, W.J.: Understanding the acceleration phenomenon via high-resolution differential equations. Math. Program. {\bf 195}, 79--148 (2022)
		
		\bibitem{bib31}
		Tikhonov, A.N.: Solution of incorrectly formulated problems and the regularization method. Sov. Dok. {\bf 4}, 1035--1038 (1963)
		
		
		
		\bibitem{bib10}
		Vassilis, A., Aujol, J.-F., Dossal, C.: The differential inclusion modeling FISTA algorithm and optimality of convergence rate in the case  $b \leq 3$. SIAM J. Optim. {\bf 28}(1), 551--574 (2018)
		
		\bibitem{bib44}
		Wibisono, A., Wilson, A.C., Jordan, M.I.: A variational perspective on accelerated methods in optimization. Proc. Natl. Acad. Sci. {\bf 113}(47), E7351--E7358 (2016)
		
		\bibitem{bib45}
		Wilson, A.C., Recht, B., Jordan, M.I.: A Lyapunov analysis of accelerated methods in optimization. J. Mach. Learn. Res. {\bf 22}, 5040--5073 (2021)
		
		
		\bibitem{bib53}
		Wang, K., Jiang, C., Ng, A.K.Y., Zhu, Z.: Air and rail connectivity patterns of major city clusters in China. Transp. Res. A: Policy. Pract. {\bf 139}, 35--53 (2020)
		
		
		
		\bibitem{bib56}
		Wan, F., Fondrevelle, J., Wang, T., Duclos, A.: Two-stage multi-objective optimization for ICU bed allocation under multiple sources of uncertainty. Sci. Rep. {\bf 13}(1), 18925 (2023)
		
		\bibitem{bib22}
		Xu, B., Wen, B.: On the convergence of a class of inertial dynamical systems with Tikhonov regularization. Optim. Lett. {\bf 15}, 2025--2052 (2021)
		
		\bibitem{bib13}
		Zeng, X., Lei, J., Chen, J.: Dynamical primal-dual Nesterov accelerated method and its application to network optimization. IEEE Trans. Automat. Control {\bf 68}(3), 1760--1767 (2022)
		
		\bibitem{bib54}
		Zheng, S., Wang, K., Dong, K., Wan, Y., Fu, X.: Does the shipping alliance aggravate or alleviate container shipping market volatility. Transp. Res. A: Policy. Pract. {\bf 189}, 104231 (2024)
		
		\bibitem{bib23}
		Zhu, T.T., Hu, R., Fang, Y.P.: Tikhonov regularized second-order plus first-order primal-dual dynamical systems with asymptotically vanishing damping for linear equality constrained convex optimization problems. Optimization (2024) \url{https://doi.org/10.1080/02331934.2024.2407515}
	\end{thebibliography}
\end{document}